\documentclass[aos,MSNbibl,seceqn,nameyear,dvips]{arximspdf}
\usepackage{subenv}
\usepackage{graphicx}

\doi{10.1214/11-AOS913}
\volume{39}
\issue{5}
\pubyear{2011}
\firstpage{2658}
\lastpage{2685}

\makeatletter

\newcommand{\cas}{\stackrel{\mathrm{a.s.}}{\longrightarrow}}
\newcommand{\cd}{\stackrel{d}{\rightarrow}}
\newcommand{\cp}{\stackrel{p}{\rightarrow}}
\newcommand{\cpshort}{\stackrel{p}{\rightarrow}}
\newcommand{\R}{\mathbb{R}}
\newcommand{\bbeta}{\bolds{\beta}}
\newcommand{\hatbbeta}{\hat{\bolds{\beta}}}
\newcommand{\vectornorm}[1]{\|#1\|}
\renewcommand{\Pr}{P}
\newcommand{\Cov}{\operatorname{Cov}}
\newcommand{\Var}{\operatorname{Var}}
\newcommand{\argmin}{\mathop{\arg\min}}
\newcommand{\matbold}{\mathbf}

\newtheorem{theorem}{Theorem}

\makeatother

\begin{document}
\begin{frontmatter}

\title{Computational approaches for empirical Bayes methods and
Bayesian sensitivity analysis}
\runtitle{Computational approaches for empirical Bayes methods}

\begin{aug}
\author[A]{\fnms{Eugenia} \snm{Buta}\corref{}\ead[label=e1]{eugenia.buta@yale.edu}}
\and
\author[B]{\fnms{Hani} \snm{Doss}\thanksref{t1}\ead[label=e2]{doss@stat.ufl.edu}}
\runauthor{E. Buta and H. Doss}
\affiliation{Yale University and University of Florida}
\address[A]{Department of Epidemiology and Public Health\\
Yale University \\
New Haven, Connecticut 06510 \\
USA\\
\printead{e1}}
\address[B]{Department of Statistics \\
University of Florida \\
Gainesville, Florida 32611 \\
USA\\
\printead{e2}} 
\end{aug}

\thankstext{t1}{Supported by NSF Grant DMS-08-05860.}

\received{\smonth{5} \syear{2010}}
\revised{\smonth{2} \syear{2011}}

%
\begin{abstract}
We consider situations in Bayesian analysis where we have a family
of priors $\nu_h$ on the parameter $\theta$, where $h$ varies
continuously over a space $\mathcal{H}$, and we deal with two related
problems. The first involves sensitivity analysis and is stated
as follows. Suppose we fix a function~$f$ of~$\theta$. How do we
efficiently estimate the posterior expectation of~$f(\theta)$
simultaneously for all~$h$ in~$\mathcal{H}$? The second problem is
how do we identify subsets of~$\mathcal{H}$ which give rise to
reasonable choices of~$\nu_h$? We assume that we are able to
generate Markov chain samples from the posterior for a finite
number of the priors, and we develop a methodology, based on a
combination of importance sampling and the use of control
variates, for dealing with these two problems. The methodology
applies very generally, and we show how it applies in particular
to a commonly used model for variable selection in Bayesian linear
regression, and give an illustration on the US crime data of
Vandaele.
\end{abstract}

%

\end{frontmatter}

\section{Introduction}
\label{secintro}

In the Bayesian paradigm we have a data vector $Y$ with density
$p_{\theta}$ for some unknown $\theta\in\Theta$, and we wish to
put a prior density on~$\theta$. The available family of prior
densities is $\{ \nu_h, h \in\mathcal{H} \}$, where $h$ is called a~%
hyperparameter. Typically, the hyperparameter is multivariate and
choosing it can be difficult. But this choice is very important and
can have a large impact on subsequent inference. There are two
issues we wish to consider:
\begin{longlist}[(A)]
\item[(A)] Suppose we fix a quantity of interest, say, $f(\theta)$,
where $f$ is a function. How do we assess how the posterior
expectation of $f(\theta)$ changes as we vary~$h$? More
generally, how do we assess changes in the posterior distribution
of~$f(\theta)$ as we vary $h$?
\item[(B)] How do we determine if a given subset of $\mathcal{H}$
constitutes a class of reasonable choices?\vadjust{\goodbreak}
\end{longlist}
The first issue is one of sensitivity analysis and the second is one
of model selection.

As an example of the kind of problem we wish to deal with, consider
the problem of variable selection in Bayesian linear regression.
Here, we have a response variable $Y$ and a set of predictors $X_1,
\ldots, X_q$, each a vector of length $m$. For every subset
$\gamma$ of $\{ 1, \ldots, q \}$ we have a potential model $\mathcal
{M}_{\gamma}$ given by
\[
Y = 1_m \beta_0 + X_{\gamma} \beta_{\gamma} + \varepsilon,
\]
where $1_m$ is the vector of $m$ $1$'s, $X_{\gamma}$ is the design
matrix whose columns consist of the predictor vectors corresponding
to the subset $\gamma$, $\beta_{\gamma}$ is the vector of
coefficients for that subset, and $\varepsilon\sim\mathcal{N}_m(0,
\sigma^2 I)$. Let $q_{\gamma}$ denote the number of variables in
the subset $\gamma$. The unknown parameter is $\theta= (\gamma,
\sigma, \beta_0, \beta_{\gamma})$, which includes the indicator of
the subset of variables that go into the linear model. A very
commonly used prior distribution on $\theta$ is given by a hierarchy
in which we first choose the indicator $\gamma$ from the
``independence Bernoulli prior''---each variable goes into the model
with a~certain probability $w$, independently of all the other
variables---and then choose the vector of regression coefficients
corresponding to the selected variables. In more detail, the model
is described as follows:
%
%
\begin{subequation}\qquad
\label{vsblm}
\begin{eqnarray}
\label{vsblm-a}
Y &\sim&\mathcal{N}_{m}(1_m \beta_0 + X_{\gamma} \beta_{\gamma},
\sigma^2 I), \\
\label{vsblm-b}
(\sigma^2, \beta_0) &\sim& p(\sigma^2, \beta_0) \propto1 /
\sigma^2; \nonumber\\[-8pt]\\[-8pt]
&&\eqntext{\mbox{given } \sigma, \beta_{\gamma} \sim
\mathcal{N}_{q_{\gamma}} ( 0, g \sigma^2 (X'_{\gamma}
X_{\gamma})^{-1} ),}\\
\label{vsblm-c}
\gamma&\sim& w^{q_{\gamma}} (1 - w)^{q - q_{\gamma}}.
\end{eqnarray}
\end{subequation}
The prior on $(\sigma, \beta_0, \beta_{\gamma})$ is Zellner's
$g$-prior introduced in \citet{Zellner1986}, and is indexed by a
hyperparameter $g$. Although this prior is improper, the resulting
posterior distribution is proper.

Note that we have used the word ``model'' in two different ways: (i)~a~model
is a specification of the hyperparameter $h$, and (ii) a
model in regression is a list of variables to include. The meaning
of the word will always be clear from context.

To summarize, the prior on the parameter $\theta= (\gamma, \sigma,
\beta_0, \beta_{\gamma})$ is given by the two-level
hierarchy (\ref{vsblm-c}) and (\ref{vsblm-b}), and is indexed by $h
= (w, g)$. Loosely speaking, when $w$ is large and $g$ is small,
the prior encourages models with many variables and small
coefficients, whereas when $w$ is small and $g$ is large, the prior
concentrates its mass on parsimonious models with large
coefficients. Therefore, the hyperparameter $h = (w, g)$ plays a
very important role, and in effect determines the model that will be
used to carry out variable selection.

A standard method for approaching model selection involves the use
of Bayes factors. For each $h \in\mathcal{H}$, let $m_h(y)$ denote
the marginal likelihood of the data under the prior $\nu_h$, that
is, $m_h(y) = \int p_{\theta}(y) \nu_h(\theta) \,d\theta$. We will\vadjust{\goodbreak}
write~$m_h$ instead of $m_h(y)$. The Bayes factor of the model
indexed by $h_2$ vs. the model indexed by $h_1$ is defined as the
ratio of the marginal likelihoods of the data under the two models,
$m_{h_2} / m_{h_1}$, and is denoted throughout by~$B(h_2, h_1)$.
Bayes factors are widely used as a criterion for comparing models in
Bayesian analyses. For selecting models that are better than others
from the family of models indexed by $h \in\mathcal{H}$, our strategy
will be to compute and subsequently compare all the Bayes factors
$B(h, h_1)$, for all $h \in\mathcal{H}$, and a~fixed hyperparameter
value $h_1$. We could then consider as good candidate models those
with values of $h$ that result in the largest Bayes factors.

Suppose now that we fix a particular function $f$ of the parameter
$\theta$; for instance, in the example, this might be the indicator
that variable $1$ is included in the regression model. It is of
general interest to determine the posterior expectation
$E_h(f(\theta) \mid Y)$ as a function of $h$ and to determine
whether or not~$E_h(f(\theta) \mid Y)$ is very sensitive to the
value of $h$. If it is not, then two individuals using two
different hyperparameters will reach approximately the same
conclusions and the analysis will not be controversial. On the
other hand, if for a function of interest the posterior expectation
varies considerably as we change the hyperparameter, then we will
want to know which aspects of the hyperparameter (e.g., which
components of $h$) produce big changes and we may want to see a plot
of the posterior expectations as we vary those aspects of the
hyperparameter. Except for extremely simple cases, posterior
expectations cannot be obtained in closed form, and are typically
estimated via Markov chain Monte Carlo (MCMC). It is slow and
inefficient to run Markov chains for every hyperparameter value $h$.
Section \ref{secbf-pe} reviews an existing method for estimating
$E_h(f(\theta) \mid Y)$ that bypasses the need to run a~separate
Markov chain for every $h$. The method has an analogue for the
problem of estimating Bayes factors. Unfortunately, the method has
severe limitations, which we also discuss.

In this paper we address the sensitivity analysis and model
selection issues discussed above. Our approach involves running
Markov chains corresponding to a few values of the hyperparameter,
say, $h_1, \ldots, h_k$, and using these to estimate $E_h(f(\theta)
\mid Y)$ for all $h \in\mathcal{H}$ and also the Bayes factors $B(h,
h_1)$ for all $h \in\mathcal{H}$. The difficulty we face is that
there is a severe computational burden caused by the requirement
that we handle a very large number of values of $h$. Our approach
for estimating large families of posterior expectations and Bayes
factors is based on a combination of MCMC, importance sampling, and
the use of control variates. The main contribution of this work is
the development of theory to support the method. This theory can be
used when dealing with implementation issues. The paper is
organized as follows. In Section \ref{secbf-pe} we describe our
methodology for estimating Bayes factors and posterior expectations,
and give statements of theoretical results associated with the
methodology. In Section \ref{secest-var-sel-sp} we discuss
estimation of the variance and implementation issues. In
Section \ref{secillustration} we return to the problem of variable
selection in Bayesian linear regression, and show how our\vadjust{\goodbreak}
methodology applies in that model. The \hyperref[app]{Appendix} gives
proofs of the
theorems stated in the paper.

The idea of doing importance sampling using data streams from
multiple densities has been investigated in several papers before.
In \citet{Vardi1985}, \citet{GillVardiWellner1988},
\citet{Geyer1994}, \citet{MengWong1996}, \citet{KongEtal2003}
and \citet{Tan2004}, it is assumed that we have samples from each
density and that each density is known except for a normalizing
constant. The objective is to estimate all possible ratios of
normalizing constants, and expectations of a given function with
respect to each of the densities. The estimates in all these papers
are identical, although the computational schemes to obtain them
given in these papers are different. \citet{GillVardiWellner1988}
and \citet{Tan2004} obtain the asymptotic distribution of the
estimates when the samples are i.i.d., and \citet{Geyer1994} gives the
asymptotic distribution when the samples are Markov chains
satisfying certain regularity conditions.

Our Bayesian framework is the same as the framework described above.
Let $\nu_{h,y}$ denote the posterior density of $\theta$ given $Y =
y$ when the prior is $\nu_h$. The posterior densities $\nu_{h_j,y}$
are given by $\nu_{h_j,y}(\theta) = p_{\theta}(y) \nu_{h_j}(\theta)
/ m_{h_j}$, where the functional form $p_{\theta}(y)
\nu_{h_j}(\theta)$ is known, but the normalizing constant $m_{h_j}$
is not. Our perspective is different from that of the previous
authors in that we are interested in estimation of the ratios
$m_h/m_{h_1}$ and of posterior expectations $\int f(\theta)
\nu_{h,y}(\theta) \,d\theta$ for a very large number of $h$'s.
Consequently, in addition to the obvious computational demands for
handling many $h$'s, we also have to deal with the fact that we will
not have a sample from $\nu_{h,y}$ for every $h \in\mathcal{H}$, but
only from $\nu_{h_j,y}, j = 1, \ldots, k$. Thus, we are
concerned with computational efficiency, in addition to statistical
efficiency. These issues are discussed in detail in
Section \ref{secbf-pe}.

\section{Estimation of Bayes factors and posterior expectations}
\label{secbf-pe}

Suppose that we have a sample $\theta_1, \ldots,
\theta_n$ (i.i.d. or ergodic Markov chain output) from the posterior
density $\nu_{h_1,y}$ for a fixed $h_1$ and we are interested in the
posterior expectation
%
%
\begin{equation}
\label{rnd}
E_h\bigl(f(\theta) \mid Y = y\bigr) = \int f(\theta)
\frac{\nu_{h,y}(\theta)}{\nu_{h_1,y}(\theta)} \nu_{h_1,y}(\theta)
\,d\theta
\end{equation}
for different values of $h$.
%
Using the fact that
\[
\int \frac{p_{\theta}(y) \nu_h(\theta)/m_h} {p_{\theta}(y)
\nu_{h_1}(\theta)/m_{h_1}} \nu_{h_1,y}(\theta) \, d\theta = 1,
\]
we see that this expectation may be written as
%
%
\begin{equation}
\label{imp-samp-deriv}
\int f(\theta) \frac{p_{\theta}(y) \nu_h(\theta) / m_h}
{p_{\theta}(y) \nu_{h_1}(\theta) / m_{h_1}} \nu_{h_1,y}(\theta)
\,d\theta= \frac{ \int f(\theta) ({\nu_h(\theta)}/
{\nu_{h_1}(\theta)}) \nu_{h_1,y}(\theta) \,d\theta} {\int
({\nu_h(\theta)}/{\nu_{h_1}(\theta)}) \nu_{h_1,y}(\theta)
\,d\theta},\hspace*{-32pt}
\end{equation}
where the right-hand side of (\ref{imp-samp-deriv}) does not involve the
ratio $m_h / m_{h_1}$. The idea to express $\int f(\theta)
\nu_{h,y}(\theta) \,d\theta$ in this way was proposed in a
different\vadjust{\goodbreak} context by \citet{Hastings1970}. The right-hand side
of (\ref{imp-samp-deriv}) is the ratio of two integrals with respect
to $\nu_{h_1,y}$, each of which may be estimated from the sequence
$\theta_1, \ldots, \theta_n$. We may estimate the numerator and the
denominator~by
%
%
\begin{equation}
\label{ratio-est}
\frac{1}{n} \sum_{i=1}^n f(\theta_i) [\nu_h(\theta_i) /
\nu_{h_1}(\theta_i)] \quad\mbox{and}\quad \frac{1}{n}
\sum_{i=1}^n [\nu_h(\theta_i) / \nu_{h_1}(\theta_i)],
\end{equation}
respectively, and $\int f(\theta) \nu_{h,y}(\theta) \,d\theta$ is
estimated by the ratio of these two quantities.

The disappearance of the likelihood function on the right-hand side
of~(\ref{imp-samp-deriv}) is very convenient because its computation
requires considerable effort in some cases (e.g., when we
have missing or censored data, the likelihood is a~possibly
high-dimensional integral). Note that the second average
in~(\ref{ratio-est}) is an estimate of $m_h / m_{h_1}$, that is, the
Bayes factor $B(h, h_1)$. Ideally, we would like to use the
estimates in (\ref{ratio-est}) for multiple values of $h$ using only
a sample from the posterior distribution corresponding to the fixed
hyperparameter value $h_1$. But, when the prior $\nu_h$ differs
from $\nu_{h_1}$ greatly, the two estimates in (\ref{ratio-est}) are
unstable because of the potential that only a~few observations will
dominate the sums. Their ratio suffers the same defect.

A natural approach for dealing with the instability of these simple
estimates is to choose $k$ values $h_1, \ldots, h_k \in\mathcal{H}$
and in (\ref{rnd}) replace $\nu_{h_1,y}$ with a~mixture
$\sum_{s=1}^k a_s \nu_{h_s,y}$, where $a_s \geq0$, for $s = 1,
\ldots, k$, and $\sum_{s=1}^k a_s = 1$. For concreteness, consider
the estimate of the Bayes factor. Let $\overline{\nu}_{\cdot y} =
\sum_{s=1}^k a_s \nu_{h_s,y}$, and let $d_s = m_{h_s} / m_{h_1},
s = 1, \ldots, k$. Note that
%
%
\begin{equation}
\label{identity-bf}
B(h, h_1) = \int\frac{\nu_h(\theta)} {\sum_{s=1}^k a_s \nu_{h_s}
(\theta)/d_s} \overline{\nu}_{\cdot y}(\theta) \,d\theta
\end{equation}
and
%
%
\begin{eqnarray}
\label{identity-pe}\quad
\int f(\theta) \nu_{h,y}(\theta)
\,d\theta& = &(B(h, h_1))^{-1} \int f(\theta) \frac{\nu_h(\theta)}
{\sum_{s=1}^k a_s \nu_{h_s} (\theta)/d_s}
\overline{\nu}_{\cdot y}(\theta) \,d\theta\nonumber\\[-8pt]\\[-8pt]
& = &\frac{ \int f(\theta)
({\nu_h(\theta)}/{\sum_{s=1}^k a_s \nu_{h_s}
(\theta)/d_s}) \overline{\nu}_{\cdot y}(\theta) \,d\theta}
{ \int({\nu_h(\theta)}/
{\sum_{s=1}^k a_s \nu_{h_s} (\theta)/d_s})
\overline{\nu}_{\cdot y}(\theta) \,d\theta}.\nonumber
\end{eqnarray}
[These two identities are valid under the condition that
$\nu_h(\theta) = 0$ whenever $\nu_{h_s}(\theta) = 0$ for all $s$.]
Suppose that\vspace*{-1pt} for each $l = 1, \ldots, k$ we have Markov chain
samples $\theta_i^{(l)}, i = 1, \ldots, n_l$, from the posterior
density $\nu_{h_l,y}$. Letting $n = \sum_{s=1}^k n_s$, if $a_s = n_s/n$, then the pooled sample is
a stratified sample from $\overline{\nu}_{\cdot y}$. \citet{Doss2010}
considers the case where the vector $d = (d_2, \ldots, d_k)'$ is
known. In this situation, the right-hand side of (\ref{identity-bf}) is
the integral of a known function with respect to the mixture density
$\overline{\nu}_{\cdot y}$. He shows that under certain regularity
conditions, the estimate of $B(h, h_1)$ obtained by replacing the
right-hand side of (\ref{identity-bf}) by its natural Monte Carlo
estimate using the pooled sample is consistent and asymptotically
normal.

In virtually all applications, the value of the vector $d$ is
unknown. The estimates of $B(h, h_1)$ and $\int f(\theta)
\nu_{h,y}(\theta) \,d\theta$ that we\vspace*{1pt} consider in this paper are
constructed by first forming an estimate $\hat{d}$ of $d$, and then
using the natural Monte Carlo estimates of the integral
in~(\ref{identity-bf}) and of the two integrals
in~(\ref{identity-pe}) with $\hat{d}$ substituted for $d$. The MCMC
scheme we will use involves the following two stages:

\begin{longlist}[\textit{Stage} 2.]
\item[\textit{Stage} 1.] Generate samples $\theta_i^{(l)0}, i = 1, \ldots,
N_l$, from $\nu_{h_l,y}$, the posterior density of $\theta$ given
$Y = y$, assuming that the prior is $\nu_{h_l}$, for each $l = 1,
\ldots, k$, and use these $N = \sum_{l=1}^k N_l$ observations to
form an estimate of $d$.
\item[\textit{Stage} 2.] Independently of stage $1$, again generate samples
$\theta_i^{(l)}, i = 1, \ldots, n_l$, from $\nu_{h_l,y}$, for
each $l = 1, \ldots, k$, and construct the estimate of the Bayes
factor $B(h, h_1)$ based on this second set of $n = \sum_{l=1}^k
n_l$ observations and the estimate of $d$ from stage $1$.
\end{longlist}
The estimate of $d$ in stage $1$ is formed using a method introduced
by \citet{Vardi1985}, and this estimate is discussed in the
beginning of Section \ref{secbf}. From now on, for $l = 1, \ldots,
k$, we use the notation $A_l$ and $a_l$ to identify the ratios~$N_l
/ N$ and $n_l / n$, respectively.

It is natural to ask why we use two steps of sampling, instead of
estimating the vector $d$ and $B(h, h_1)$ from a single sample. The
quantity considered in \citet{Doss2010} is
%
%
\begin{equation}
\label{est-bf-simple}
\hat{B}(h, h_1, d) = \sum_{l=1}^k \sum_{i=1}^{n_l} \frac{\nu_h
(\theta_i^{(l)})} {\sum_{s=1}^k n_s \nu_{h_s} (\theta_i^{(l)}) /
d_s},
\end{equation}
and it involves the vector $d$. The estimate considered in the present
paper is $\hat{B}(h, h_1, \hat{d})$, where $\hat{d}$ is an estimate of
$d$. The variance of $\hat{B}(h, h_1, \hat{d})$ turns out to be greater
than that of $\hat{B}(h, h_1, d)$ (and this\vspace*{1pt} is true
whether we use two steps of sampling or a single step).
Thus,\vspace*{1pt} the variance decomposes as $\Var( \hat{B}(h, h_1,
\hat{d}) ) = \Var( \hat{B}(h, h_1, d) ) + V_d$, where $V_d$ is the
increase in variance resulting from using $\hat{d}$ instead of $d$.
Because we wish to estimate $B(h, h_1)$ for a large number of $h$'s and
for each $h$ the computational time needed is linear in the total
sample size, this total sample size cannot be very large. On the other
hand, $d$ needs to be estimated only once. So if generating the chains
is not computationally demanding, then one can use very long chains to
estimate $d$ and so greatly reduce the term $V_d$. A precise statement
regarding the benefits of the two-stage scheme would have to take into
account the cost of computing the typical term in (\ref{est-bf-simple})
and the cost of generating a point in the chain, and no such statement
can be made at the level of generality considered in this paper.
However, in all the examples we have encountered, for fixed
computational resources, the two-stage scheme gives estimates with
considerably smaller variance. We mention here that our theoretical
results are stated for the two-stage schemes, but these results have
analogues for the case where a single sample is used to estimate both
$d$ and the family of Bayes factors $B(h, h_1), h \in\mathcal{H}$, and
these are given in \citet{Buta2010}.

A summary of the main contributions of the present work is as
follows:

{\renewcommand\thelonglist{\arabic{longlist}}
\begin{longlist}
\item We develop a complete characterization of the asymptotic
distribution of the estimate (\ref{est-bf-simple}) and also of a
variant involving the use of control variates developed by
\citet{Doss2010} for the realistic case where $d$ is estimated
from stage $1$ sampling. Included in our results is an explicit formula
for the increase in variance resulting from using an estimate of $d$
instead of $d$ itself. (This contradicts statements in the literature
to the effect that using a $\sqrt{n}$-consistent estimate of $d$ rather
than $d$ itself does not inflate the variance; see our discussion in
the \hyperref[app]{Appendix}.)
\item We develop an analogous theory for the problem of estimating a
family of posterior expectations $E_h(f(\theta) \mid Y = y), h
\in\mathcal{H}$.
\item For any of our estimators, the variance is a sum of two
components, and we discuss how each of these may be estimated. An
important problem is how to properly select the skeleton points
$h_1, \ldots, h_k$, and ideally we would like to position these in
such a way that the variance is minimized. We show how the
variance estimates can be used to suggest good sets of skeleton
points.
\item We apply the methodology to the problem of Bayesian variable
selection discussed earlier. In particular, we show how our
methods enable us to select good values of $h = (w, g)$ and to
also see how the probability that a~given variable is included in
the regression varies with $(w, g)$.
\end{longlist}}

\subsection{Estimation of Bayes factors}
\label{secbf}
Here, we analyze the asymptotic distributional properties of the
estimator that results if in (\ref{est-bf-simple}) we replace $d$
with an estimate. \citet{Geyer1994} proposes an estimator for $d$
based on the ``reverse logistic regression'' method and Theorem $2$
therein shows that this estimator is asymptotically normal when the
samplers used satisfy certain regularity conditions. This estimator
is obtained by maximizing with respect to $d_2, \ldots, d_k$ the log
quasi-likelihood
%
%
\begin{equation}
\label{geyer-criterion}
l_N(d) = \sum_{l=1}^k \sum_{i=1}^{N_l} \log\Biggl( \frac{ A_l
\nu_{h_l} (\theta_i^{(l)0}) / d_l } {\sum_{s=1}^k A_s \nu_{h_s}
(\theta_i^{(l)0}) / d_s} \Biggr).
\end{equation}
As was mentioned earlier, the estimate is the same as the estimates
obtained by \citet{Vardi1985}, \citet{MengWong1996} and
\citet{KongEtal2003}. We assume that for all the Markov chains we
use a Strong Law of Large Numbers (SLLN) holds for all integrable
functions [for sufficient conditions see, e.g., Theorem~$2$ of
\citet{AthreyaDossSethuraman1996}]. In the next theorem we show
that if $\hat{d}$ is the estimate produced by
Geyer's (\citeyear{Geyer1994})
method, or any of the equivalent estimates discussed above, then the
estimate of the Bayes factor given by
%
%
\begin{equation}
\label{est-bf-complex}
\hat{B}(h, h_1, \hat{d}) = \sum_{l=1}^k \sum_{i=1}^{n_l} \frac{
\nu_h (\theta_i^{(l)}) } { \sum_{s=1}^k n_s \nu_{h_s}
(\theta_i^{(l)}) / \hat{{d}_{s}} }
\end{equation}
is asymptotically normal if certain regularity conditions are met.
In (\ref{est-bf-complex}), $\hat{d}_1 = 1$.

Before we state the theorem, we need to define the expressions that
appear in the asymptotic variance. For $l = 1, \ldots, k, i = 1,
\ldots, n_l$, let
%
%
\begin{equation}
\label{reg-prob-y}
Y_{i,l} = \frac{ \nu_h(\theta_i^{(l)}) } { \sum_{s=1}^k a_s
\nu_{h_s} (\theta_i^{(l)}) / d_s }
\end{equation}
(the $Y_{i,l}$'s depend on $h$, but this dependence is suppressed to
lighten the notation), and let
\[
\tau_l^2(h) = \operatorname{Var}(Y_{1,l}) + 2 \sum_{g=1}^{\infty}
\operatorname{Cov}(Y_{1,l}, Y_{1+g,l}),\qquad \tau^2(h) = \sum_{l=1}^k
a_l \tau_l^2(h).
\]
Also, let $c(h)$ be the vector of length $k - 1$ for which the $(j -
1)$th coordinate is
%
%
\begin{eqnarray}
\label{c}
[c(h)]_{j-1} = \frac{B(h, h_1)}{d_j^2} \int\frac{ a_j
\nu_{h_j}(\theta) } { \sum_{s=1}^k a_s \nu_{h_s} (\theta) / d_s }
\cdot\nu_{h,y}(\theta) \,d\theta,\nonumber\\[-8pt]\\[-8pt]
&&\eqntext{\qquad j = 2, \ldots, k.}
\end{eqnarray}
\begin{theorem}
\label{thmbf-d-est}
Let $h \in\mathcal{H}$ be fixed. Suppose the chains in stage $2$
satisfy conditions \textup{(A1)} and \textup{(A2)} in Doss
(\citeyear{Doss2010}):

{\renewcommand\thelonglist{(A\arabic{longlist})}
\renewcommand\labellonglist{\thelonglist}
\begin{longlist}
\item
\label{A1}
For each $l = 1, \ldots, k$, the chain $\{ \theta_i^{(l)}
\}_{i=1}^{\infty}$ is geometrically ergodic.
\item
\label{A2}
For each $l = 1, \ldots, k$, there exists $\varepsilon> 0$ such
that
%
%
\begin{equation}
\label{mom-cond0}
E\biggl( \biggl| \frac{ \nu_h(\theta_1^{(l)}) } {
\sum_{s=1}^k a_s \nu_{h_s}(\theta_1^{(l)}) / d_s } \biggr|^{2 +
\varepsilon} \biggr) < \infty.
\end{equation}
\end{longlist}}

\noindent In the expectation in (\ref{mom-cond0}), $\theta_1^{(l)} \sim
\nu_{h_l,y}$. Assume\vspace*{-1pt} also that the chains in stage~$1$
satisfy the conditions in\vspace*{1pt} Theorem $2$ of Geyer (\citeyear{Geyer1994})
that imply $\sqrt{N} (\hat{d} - d) \cd\mathcal{N}(0, \Sigma)$. In
addition, suppose the total sample sizes for the two stages, $N$ and
$n$, satisfy $n \rightarrow\infty$, and $N \rightarrow\infty$ in such a
way that $n / N \rightarrow q \in[0, \infty)$. Then
\[
\sqrt{n} \bigl( \hat{B}(h, h_1, \hat{d}) - B(h, h_1) \bigr) \cd
\mathcal{N}\bigl( 0, q {c(h)}' \Sigma c(h) + \tau^2(h) \bigr).
\]
\end{theorem}

As alluded to earlier, there are two components to the expression
for the variance. The first component arises from estimating $d$,
and the second component is the variance that we would have if we
had estimated the Bayes factor knowing what $d$ is. As can be seen
from the formula, the first component vanishes if $q = 0$, that is, if
the sample size for estimating the parameter~$d$ converges to
infinity at a faster rate than does the sample size used to estimate
the Bayes factor. In this case the Bayes factor
estimator~(\ref{est-bf-complex}) using the estimate $\hat{d}$ has
the same asymptotic distribution as the estimator
in~(\ref{est-bf-simple}) which uses the true value of $d$.
Otherwise, the variance of~(\ref{est-bf-complex}) is greater than
that of~(\ref{est-bf-simple}), and the difference between the
variances depends on the parameter $q$. This parameter is
determined by the user and should be chosen in such a way as to
minimize the variance given computer resources; this is discussed in
Section~\ref{secest-var-sel-sp}.

\subsection{Estimation of Bayes factors using control variates}
\label{seccv}
Recall that we have samples $\theta_i^{(l)}, i = 1, \ldots, n_l$,
from $\nu_{h_l,y}, l = 1, \ldots, k$, with independence across
samples (stage $2$ of sampling) and that, based on an independent
set of preliminary MCMC runs (stage $1$ of sampling), we have
estimated the constants $d_2, \ldots, d_k$. Also, $n_l/n = a_l$ and
$n = \sum_{l=1}^k n_l$. Let
%
%
\begin{equation}
\label{def-Y}
Y(\theta) = \frac{ \nu_h(\theta) } { \sum_{s=1}^k a_s \nu_{h_s}
(\theta) / d_s }.
\end{equation}
Recalling that $\overline{\nu}_{\cdot y} := \sum_{s=1}^k a_s
\nu_{h_s,y}$, we have $E_{\overline{\nu}_{\cdot y}} (Y(\theta)) = B(h,
h_1)$, where the subscript $\overline{\nu}_{\cdot y}$ to the expectation
indicates that $\theta\sim\overline{\nu}_{\cdot y}$. Also, for $j = 2,
\ldots, k$, let
%
%
\begin{eqnarray}
\label{def-cv-a}
Z^{(j)}
(\theta) & = & \frac{ \nu_{h_j}(\theta) / d_j - \nu_{h_1} (\theta) }
{ \sum_{s=1}^k a_s \nu_{h_s}(\theta) / d_s } \\
\label{def-cv-b}
& = & \frac{ \nu_{h_j,y}(\theta) - \nu_{h_1,y}(\theta) } {
\sum_{s=1}^k a_s \nu_{h_s,y} (\theta) }.
\end{eqnarray}
Expression (\ref{def-cv-b}) shows that $E_{\overline{\nu}_{\cdot y}}
(Z^{(j)}(\theta)) = 0$. This is true even if the priors $\nu_{h_j}$
and $\nu_{h_1}$ are improper, as long as the posteriors
$\nu_{h_j,y}$ and $\nu_{h_1,y}$ are proper, exactly our situation in
the Bayesian variable selection example of Section \ref{secintro}.
On the other hand, the representation (\ref{def-cv-a}) shows that
$Z^{(j)} (\theta)$ is computable if we know the $d_j$'s---it
involves the priors and not the posteriors. [A similar remark
applies to (\ref{def-Y}).] Therefore, if as in\vadjust{\goodbreak} \citet{Doss2010} we
define for $l = 1, \ldots, k, i = 1, \ldots, n_l$
%
%
\begin{equation}
\label{reg-prob-z}\quad
Z_{i,l}^{(1)} = 1,\qquad Z_{i,l}^{(j)} = \frac{
\nu_{h_j}(\theta_i^{(l)}) / d_j - \nu_{h_1}(\theta_i^{(l)}) } {
\sum_{s=1}^k a_s \nu_{h_s} (\theta_i^{(l)}) / d_s },\quad j = 2,
\ldots, k,
\end{equation}
then for any fixed $\bbeta= (\beta_2, \ldots, \beta_k)$,
%
%
\begin{equation}
\label{ihat-beta-d}
\hat{I}_{\bbeta}^{d} = \frac{1}{n} \sum_{l=1}^k \sum_{i=1}^{n_l}
\Biggl( Y_{i,l} - \sum_{j=2}^k \beta_{j} Z_{i,l}^{(j)}
\Biggr)
\end{equation}
is an\vspace*{1pt} unbiased estimate of $B(h, h_1)$. The value of $\bbeta$ that
minimizes the variance of $\hat{I}_{\bbeta}^{d}$ is unknown. As is
commonly done when one uses control variates, we use instead the
estimate obtained by doing ordinary linear regression of the
response $Y_{i,l}$ on the predictors $Z_{i,l}^{(j)}, j = 2,
\ldots, k$, and to emphasize that this estimate depends on~$d$, we
denote it by $\hatbbeta(d)$. \citet{Doss2010} shows that~%
$\hatbbeta(d)$ converges almost surely to a finite limit,
$\bbeta_{\lim}$. His Theorem $1$ states that the estimator
$\hat{B}_{\mathrm{reg}}(h, h_1) = \hat{I}_{\hatbbeta(d)}^{d}$,
obtained under the assumption that we know the constants $d_2,
\ldots, d_k$, has an asymptotically normal distribution. As
mentioned earlier, $d_2, \ldots, d_k$ are typically unknown, and
must be estimated. Let $\hat{d}_2, \ldots, \hat{d}_k$ be estimates
obtained from previous MCMC runs and let
%
%
\begin{equation}
\label{est-bf-complex-cv}
\hat{I}_{{\hatbbeta(\hat{d})}}^{\hat{d}} = \frac{1}{n}
\sum_{l=1}^k \sum_{i=1}^{n_l} \Biggl( \hat{Y}_{i,l} - \sum_{j=2}^k
\hat{\beta}_{j}(\hat{d}) \hat{Z}_{i,l}^{(j)} \Biggr),
\end{equation}
where $\hat{Y}_{i,l}$ and $\hat{Z}_{i,l}^{(j)}$ are as
in (\ref{reg-prob-y}) and (\ref{reg-prob-z}), except using $\hat{d}$
for $d$, and $\hatbbeta(\hat{d})$ is the least squares regression
estimator from regressing $\hat{Y}_{i,l}$ on predictors~%
$\hat{Z}_{i,l}^{(j)}$, $j = 2, \ldots, k$.

The next theorem gives the asymptotic distribution of this new
estimator, and before we state it we introduce some notation. Let
%
%
\begin{equation}
\label{Uil}
U_{i,l} = Y_{i,l} - \sum_{j=2}^k \beta_{j,\lim}
Z_{i,l}^{(j)}
\end{equation}
and let
%
%
\begin{equation}
\label{sigmah}
\sigma_l^2(h) = \operatorname{Var}(U_{1,l}) + 2 \sum_{g=1}^{\infty}
\operatorname{Cov}(U_{1,l}, U_{1+g,l}),\qquad \sigma^2(h) = \sum_{l=1}^k
a_l \sigma_l^2(h).\hspace*{-32pt}
\end{equation}
Also, let $w(h)$ be the vector of length $k - 1$ for which the $(t -
1)$th coordinate ($t = 2, \ldots, k$) is
%
%
\begin{eqnarray}
\label{grad-K-paper}
[w(h)]_{t-1} & = &\frac{B(h, h_1)} {d_t^2} \int\frac{ a_t
\nu_{h_t}(\theta) } { \sum_{s=1}^k a_s
\nu_{h_s} (\theta) / d_s } \cdot
\nu_{h,y}(\theta) \,d\theta+
\beta_{t,\lim} \frac{1}{d_t} \nonumber\\
&&{} + \sum_{j=2}^k \beta_{j,\lim}
\int\frac{ a_t \nu_{h_t}(\theta) } { d_t^2
\sum_{s=1}^k a_s \nu_{h_s}(\theta) / d_s }
\\
&&\hphantom{+ \sum_{j=2}^k \beta_{j,\lim}
\int}
{}\times\bigl( \nu_{h_1,y}(\theta) -
\nu_{h_j,y}(\theta) \bigr) \,d\theta.
\nonumber
\end{eqnarray}

\begin{theorem}
\label{thmbf-d-est-cv}
Suppose all the conditions from Theorem \ref{thmbf-d-est} are
satisfied. Moreover, assume that $\matbold{R}$, the $k \times
k$ matrix defined by
\[
R_{j,j'} = E\Biggl( \sum_{l=1}^k a_l Z_{1,l}^{(j)}
Z_{1,l}^{(j')} \Biggr),\qquad j, j' = 1, \ldots, k,
\]
is nonsingular. Then
\[
\sqrt{n} \bigl( \hat{I}_{{\hatbbeta(\hat{d})}}^{\hat{d}} - B(h,
h_1) \bigr) \cd\mathcal{N}\bigl( 0, q {w(h)}' \Sigma w(h) +
\sigma^2(h) \bigr).
\]
\end{theorem}

As mentioned above, for any $\bbeta$, $\hat{I}_{\bbeta}^{d}$
in (\ref{ihat-beta-d}) is an unbiased estimate of $B(h, h_1)$, which
leads to the question of what is the optimal value of $\bbeta$ to
use. It is not difficult to see that when each of the sequences $\{
\theta_i^{(l)} \}_{i=1}^{n_l}$ is i.i.d., the value of $\bbeta$ that
minimizes the variance of $\hat{I}_{\bbeta}^{d}$ is
\[
\bbeta_{\mathrm{opt},\mathrm{i.i.d.}} := \argmin_{\bbeta}
\Var_{\overline{\nu}_{\cdot y}} \Biggl( Y(\theta) - \sum_{j=2}^k
\beta_{j} Z^{(j)}(\theta) \Biggr),
\]
that is, the optimal value is the same whether we have a random sample
from $\overline{\nu}_{\cdot y}$ or a stratified sample. It is natural to
ask whether $\bbeta_{\mathrm{opt},\mathrm{i.i.d.}}$ is still optimal
when the
$k$ sequences $\{ \theta_i^{(l)} \}_{i=1}^{n_l}$ are Markov chains.
It turns out that:
\begin{longlist}
\item$\bbeta_{\mathrm{opt},\mathrm{i.i.d.}}$ is not optimal,
\item using $\bbeta_{\mathrm{opt},\mathrm{i.i.d.}}$ can actually
increase the
variance (when the Markov chains mix at significantly different
rates, chains that are of the same length do not have the same
``effective sample sizes,'' but $\bbeta_{\mathrm{opt},\mathrm{i.i.d.}}$ does
not reflect this fact).
\end{longlist}
In our experience, using $\bbeta_{\mathrm{opt},\mathrm{i.i.d.}}$, or, more
precisely, the least squares estimate [which in \citet{Doss2010} was
shown to converge almost surely to $\bbeta_{\mathrm{opt},\mathrm{i.i.d.}}$],
typically gives a significant reduction in variance.
\citet{ButaDoss2011} prove points (i) and (ii) above and also
discuss an approach for estimating the value of $\bbeta$ that is
optimal in the Markov chain case.

\subsection{Estimation of posterior expectations}
\label{secpe}
In this section we describe a~method for estimating the posterior
expectation of a function $f$ when the prior is $\nu_h$. Let us
denote this quantity by
\[
I^{[f]}(h) = \int f(\theta)\nu_{h,y}(\theta) \,d\theta.
\]
Define
\begin{eqnarray*}
Y_{i,l}^{[f]}& = &\frac{ f(\theta_i^{(l)}) \nu_h (\theta_i^{(l)}) }
{ \sum_{s=1}^k a_s \nu_{h_s} (\theta_i^{(l)}) /
d_s } = \frac{ f(\theta_i^{(l)}) \nu_h (\theta_i^{(l)}) / m_h }
{ \sum_{s=1}^k a_s \nu_{h_s} (\theta_i^{(l)}) / m_{h_s}
} \cdot\frac{m_h}{m_{h_1}} \\[-2pt]
& = & \frac{ f(\theta_i^{(l)}) \nu_{h,y} (\theta_i^{(l)}) } {
\sum_{s=1}^k a_s \nu_{h_s,y} (\theta_i^{(l)}) } B(h,
h_1).
\end{eqnarray*}
With the\vspace*{1pt} view of applying identity (\ref{identity-pe}), we note
that, assuming a SLLN holds for the Markov chains $\theta_i^{(l)}$,
$l = 1, \ldots, k, i = 1, \ldots, n_l$, we have
\begin{eqnarray*}
\frac{1}{n} \sum_{l=1}^k \sum_{i=1}^{n_l}
Y_{i,l}^{[f]} &=& \sum_{l=1}^k \frac{1}{n_l} \sum_{i=1}^{n_l}
\frac{n_l}{n} Y_{i,l}^{[f]} \\[-2pt]
&\cas&\int\frac{ f(\theta)\nu_{h,y} (\theta) } {
\sum_{s=1}^k a_s \nu_{h_s,y} (\theta) }
\sum_{l=1}^k a_l \nu_{h_l,y} (\theta) \,d\theta
\cdot B(h, h_1) \\[-2pt]
&=& I^{[f]}(h) \cdot B(h, h_1)
\end{eqnarray*}
and
\[
\frac{1}{n} \sum_{l=1}^k \sum_{i=1}^{n_l} Y_{i,l} \cas B(h, h_1).
\]
[The $Y_{i,l}$'s are defined in (\ref{reg-prob-y}); note that
$Y_{i,l} = Y_{i,l}^{[f]}$ when $f \equiv1$.] Letting
%
%
\begin{equation}
\label{Ihat-f-hd}
\hat{I}^{[f]}(h, d) = \frac{ \sum_{l=1}^k \sum_{i=1}^{n_l}
Y_{i,l}^{[f]} } { \sum_{l=1}^k \sum_{i=1}^{n_l} Y_{i,l} },
\end{equation}
we see that $\hat{I}^{[f]}(h, d) \cas I^{[f]}(h)$, and replacing $d$
with the estimate $\hat{d}$ obtained from stage $1$ sampling, we
form
%
%
\begin{equation}
\label{Ihat-f-h-d-est}
\hat{I}^{[f]}(h, \hat{d}) = \frac{ \sum_{l=1}^k
\sum_{i=1}^{n_l} { f(\theta_i^{(l)}) \nu_h (\theta_i^{(l)}) }
/({ \sum_{s=1}^k a_s \nu_{h_s} (\theta_i^{(l)}) / \hat{d}_{s} }) } {
\sum_{l=1}^k \sum_{i=1}^{n_l} { \nu_h
(\theta_i^{(l)}) }/ ({ \sum_{s=1}^k a_s \nu_{h_s} (\theta_i^{(l)}) /
\hat{d}_{s} }) }.
\end{equation}

The following theorem concerns the asymptotic behavior of this
estimator, and to state it, we first define the expressions that
appear in the asymptotic variance. Let
\begin{eqnarray*}
\gamma_{11} & = & \Var\bigl( Y_{1,l}^{[f]} \bigr) + 2
\sum_{g=1}^{\infty} \Cov\bigl( Y_{1,l}^{[f]},
Y_{1+g,l}^{[f]} \bigr), \\[-2pt]
\gamma_{12} & = & \gamma_{21}  = \Cov\bigl( Y_{1,l}^{[f]}, Y_{1,l}
\bigr) + \sum_{g=1}^{\infty} \bigl[
\Cov\bigl( Y_{1,l}^{[f]}, Y_{1+g,l} \bigr) +
\Cov\bigl( Y_{1,l}, Y_{1+g,l}^{[f]} \bigr)
\bigr], \\[-2pt]
\gamma_{22} & = & \Var( Y_{1,l} ) + 2
\sum_{g=1}^{\infty} \Cov( Y_{1,l},
Y_{1+g,l} )
\end{eqnarray*}
%
and
%
%
\begin{equation}
\label{Gamma}
\Gamma_l(h) =
\pmatrix{
\gamma_{11} & \gamma_{12} \cr
\gamma_{21} & \gamma_{22}}
,\qquad \Gamma(h) = \sum_{l=1}^k a_l \Gamma_l(h).
\end{equation}
Since (\ref{Ihat-f-hd}) and (\ref{Ihat-f-h-d-est}) are ratios to
which we will apply the delta method, we will consider the function
$g(u, v) = u/v$, whose gradient is $\nabla g(u, v) = (1/v, -u /
v^2)'$. Let
%
%
\begin{eqnarray}
\label{rho}
\rho(h) &=& \nabla g\bigl( I^{[f]}(h)B(h, h_1), B(h, h_1)
\bigr)'\nonumber\\[-8pt]\\[-8pt]
&&{}\times\Gamma(h) \cdot\nabla g\bigl( I^{[f]}(h) B(h, h_1), B(h,
h_1) \bigr).\nonumber
\end{eqnarray}
Finally, let $v(h)$ be the vector of length $k - 1$ for which the
$(j - 1)$th coordinate is
%
%
\begin{eqnarray}
\label{cas-grad-pe-paper}
[v(h)]_{j-1} = \int\frac{ [f(\theta) - I^{[f]}(h)] a_j
\nu_{h_j}(\theta) /d_j^2 } { \sum_{s=1}^k a_s \nu_{h_s} (\theta) /
d_s } \nu_{h,y}(\theta) \,d\theta,\nonumber\\[-8pt]\\[-8pt]
&&\eqntext{j = 2, \ldots, k.}
\end{eqnarray}
\begin{theorem}
\label{thmpe-d-est}
Suppose the conditions stated in Theorem \ref{thmbf-d-est} are
satisfied and, in addition, for each $l = 1, \ldots, k$, there
exists an $\varepsilon> 0$ such that
%
%
\begin{equation}
\label{mom-cond}
E\bigl( \bigl| Y_{1,l}^{[f]} \bigr|^{2 + \varepsilon} \bigr) <
\infty.
\end{equation}
Then
\[
\sqrt{n} \bigl( \hat{I}^{[f]} (h, \hat{d}) - I^{[f]}(h) \bigr)
\cd\mathcal{N} \bigl( 0, q v(h)' \Sigma v(h) + \rho(h) \bigr).
\]
\end{theorem}

The numerator of $\hat{I}^{[f]} (h, \hat{d})$ is an estimate of
$I^{[f]}(h) B(h, h_1)$ and the denominator is an estimate of
$B(h, h_1)$. It is possible to adjust both the numerator and
denominator through the use of control variates and thus arrive at a
variant of $\hat{I}^{[f]} (h, \hat{d})$; the theory for this is
developed in \citet{Buta2010}. As for the case of estimating the
Bayes factors, the variant is not guaranteed to give an improvement,
but a large improvement is often noted.

\section{Variance estimation and selection of the skeleton
points}
\label{secest-var-sel-sp}
Estimation of the variance of our estimates is important for several
reasons. In addition to the usual need for providing error margins
for our point estimates, variance estimates are of great help in
selecting the skeleton points. The main approaches for estimation
of the variance are (i) spectral methods, (ii)~methods based on
batching, and (iii) methods based on regeneration; see
\citet{FlegalJones2010} and \citet{MyklandTierneyYu1995} for a
review. Methods based on batching are difficult to use in our
framework because of two complications, namely, that we are dealing
with multiple chains, and we have a two-stage scheme; and procedures
based on regeneration are often difficult to implement. Here we
describe a way of estimating the variance using spectral methods.

For the sake of concreteness, consider $\hat{B}(h, h_1, \hat{d})$,
whose asymptotic variance is the expression $\kappa^2(h) = q {c(h)}'
\Sigma c(h) + \tau^2(h)$ (see Theorem \ref{thmbf-d-est}). The
term\vspace*{1pt}
$\tau^2(h)$ is the asymptotic variance of the quantity $\hat{B}(h,
h_1, d)$ in (\ref{est-bf-simple}), and since the $k$ Markov chains
are independent, $\tau^2(h) = \sum_{l=1}^k a_l \tau^2_l(h)$, whe\-re~%
$\tau^2_l(h)$ is the asymptotic variance of
%
%
\begin{equation}
\label{avg-d}
\frac{1}{n_l} \sum_{i=1}^{n_l} \frac{ \nu_h(\theta_i^{(l)}) } {
\sum_{s=1}^k a_s \nu_{h_s} (\theta_i^{(l)}) / d_s }.
\end{equation}
Now for each $l$ we will estimate $\tau^2_l(h)$ by the asymptotic
variance of
%
%
\begin{equation}
\label{avg-dhat}
\frac{1}{n_l} \sum_{i=1}^{n_l} \frac{ \nu_h(\theta_i^{(l)}) } {
\sum_{s=1}^k a_s \nu_{h_s} (\theta_i^{(l)}) / \hat{d}_s },
\end{equation}
where $\hat{d}$ is formed from stage $1$ runs. It is not too
difficult to show that under our asymptotic regime where $n / N
\rightarrow q \in[0, \infty)$, standard consistent spectral
estimates of the asymptotic variance of (\ref{avg-dhat}) are also
consistent estimates of the asymptotic variance of (\ref{avg-d});
details are given in \citet{ButaDoss2011}. \citet{Geyer1994}
gives an expression for $\Sigma$ that is explicit enough to enable
us to estimate it via standard spectral methods. Now, $c(h)$ is a~%
vector each of whose components is an integral with respect to the
posterior~$\nu_{h,y}$ [see~(\ref{c})]. The estimate derived in
Section~\ref{secpe} [see~(\ref{Ihat-f-h-d-est})] is designed
precisely to estimate such posterior expectations. Combining, we
arrive at an overall estimate of $\kappa^2(h)$, and the asymptotic
variances of our other estimates are handled similarly.

\subsection*{Selection of the skeleton points}
The asymptotic variances of any of our estimates depend\vspace*{1pt} on
the choice of the points $h_1, \ldots, h_k$. For concreteness, consider
$\hat{B}(h, h_1, \hat{d})$, and to emphasize this\vspace*{1pt}
dependence, let $V(h, h_1, \ldots, h_k)$ denote the asymptotic variance
of $\hat{B}(h, h_1, \hat{d})$. For fixed $h_1, \ldots, h_k$,
identifying the set of $h$'s for which $V(h, h_1, \ldots, h_k)$ is
\textit{finite} is typically a feasible problem. For instance,
\citet{DossTierney-disc1994} considered the pump data example
discussed in \citet{Tierney1994}, for which the hyperparameter $h$
has dimension $3$, and determined this set for the case $k = 1$. He
showed that one can go as far away from $h_1$ as one wants in certain
directions, but in other directions the range is limited. (The
calculation can be extended to any $k$.) Suppose now that we fix a
range $\mathcal{H}$ over which $h$ is to vary. A necessary first step
is to select $h_1, \ldots, h_k$ such that $V(h, h_1, \ldots, h_k) <
\infty$ for all $h \in\mathcal{H}$. Typically, however, we will want
more, and we will face the problem below.

\textit{Design problem}: find the values of the skeleton points
$h_1, \ldots, h_k$ that minimize $\max_{h \in\mathcal{H}} V(h, h_1,
\ldots, h_k)$.

Unfortunately, except for extremely simple cases, it is not possible
to calculate $V(h, h_1, \ldots, h_k)$ analytically [even if $k = 1$,
$V(h, h_1)$ is an infinite sum each of whose terms depends on the
Markov transition function in a complicated way], and maximizing
it over $h \in\mathcal{H}$ would present additional difficulties.
Furthermore, even if we were able to calculate $\max_{h \in\mathcal
{H}} V(h, h_1, \ldots, h_k)$, the design problem would involve the
minimization of a function of $k \times\dim(\mathcal{H})$ variables,
and, in general, solving the design problem is hopeless.

In our experience, we have found that the following method works
reasonably well. Having specified the range $\mathcal{H}$, we select
trial values $h_1, \ldots, h_k$ and plot the estimated variance as a
function of $h$, using one of the methods described above. If we
find a region in $\mathcal{H}$ where this variance is unacceptably
large, we ``cover'' this region by moving some $h_l$'s closer to the
region, or by simply adding new $h_l$'s in that region, which
increases $k$. This is illustrated in the example in
Section \ref{secillustration}.

\subsection*{The relative lengths of the stages 1 and 2 chains}
The parameter $q$ affects the performance of any of the methods, and
the optimal value involves a~trade-off between time spent calculating
density ratios in stage $2$ and time spent generating the chains in
stage $1$. Consider, for instance, the
estimate~(\ref{est-bf-complex}), whose asymptotic variance is given
by Theorem \ref{thmbf-d-est} and which we will write as
$\kappa^2(h) = q v_1(h) + v_2(h)$. In the discussion below, we
assume that we have run a small pilot experiment that has enabled us
to adequately estimate the components $v_1(h)$ and $v_2(h)$, and we
assume that the total sample sizes~$n$ and $N$ are both large. The
discussion is heuristic in that we assume that~$v_1(h)$ and $v_2(h)$
are nearly constant in $h$. Let $t_1$ denote the time it typically
takes to generate a single step in a chain, let $t_2$ denote the
time it takes to compute the typical term in (\ref{est-bf-complex}),
and let $g$ denote the number of values in~$\mathcal{H}$ for which we
wish to compute the estimate (\ref{est-bf-complex}). Suppose we are
given a computational budget of $T$ units of time. For any $q \in
(0, \infty)$, the time it takes to compute (\ref{est-bf-complex})
for $g$ values of $h$ is $t(q) = (n/q) t_1 + n t_1 + ng t_2$, and
setting this equal to $T$ determines $n$ to be $qT / ( (q +
1)t_1 + qg t_2 )$. The variance of the estimate is then $V(q)
= T^{-1} ( v_1(h) + v_2(h) / q ) ( (q + 1)t_1 + qg
t_2 )$. Clearly, $V(q)$~is unbounded as $q \rightarrow0$ or $q
\rightarrow\infty$. The function has a~unique minimum, which
occurs at $q_{\mathrm{opt}} = \sqrt{ [v_2(h) t_1] / [{v_1(h) (t_1 +
gt_2)}]}$. This last formula expresses in a usable manner the
intuitive notion that if $g$ is large, or if the cost of evaluating
the density ratios in (\ref{est-bf-complex}) is high relative to the
cost of running the chains, then a small value of $q$ should be
used.

\section{Illustration on variable selection in Bayesian linear
regression}
\label{secillustration}
There exist many classes of problems in Bayesian analysis in which
the sensitivity analysis and model selection issues discussed
earlier arise; see Section \ref{secdisc}. Here we give an
illustration involving the hierarchical prior used in variable
selection in the Bayesian linear regression model discussed in
Section \ref{secintro}. For this model, the parameter is the
vector $\theta= (\gamma, \sigma, \beta_0, \beta_{\gamma})$, and the
prior on~$\theta$ is given by the hierarchy (\ref{vsblm-c})
and (\ref{vsblm-b}). There exist several MCMC-based methods for
estimating the posterior distribution of $\theta$ given $Y = y$, and
the algorithm we use here is based on the Gibbs sampler of
\citet{SmithKohn1996}, which runs on the space of model indicators.
Our algorithm, developed in \citet{Buta2010}, is a Markov chain on
$\theta$ that is uniformly ergodic and also computationally
efficient (it avoids the need for repeated time-consuming matrix
inversion). It is implemented in the R package \texttt{bvslr},
available from \url{http://www.stat.ufl.edu/\textasciitilde ebuta/BVSLR}.

In Sections \ref{secintro} and \ref{secbf-pe}, $\nu_h$ and
$\nu_{h,y}$ refer to the prior and posterior \textit{densities}, and
all estimates in Section \ref{secbf-pe} involve ratios of these
prior densities. In the Bayesian linear regression model that we
are considering here, the priors $\nu_h$ on $(\gamma, \sigma,
\beta_0, \beta_{\gamma})$ are actually probability measures on $\{
0, 1 \}^q \times(0, \infty) \times\R^{q+1}$, which in fact are not
absolutely continuous with respect to the product\vspace*{0.5pt} of counting
measure on $\{ 0, 1 \}^q$ and Lebesgue measure on $(0, \infty)
\times\R^{q+1}$. For $h_1 = (w_1, g_1)$ and $h_2 = (w_2, g_2)$,
the Radon--Nikodym derivative of $\nu_{h_1}$ with respect to
$\nu_{h_2}$ is given by
%
%
\begin{eqnarray}
\label{lr}
\biggl[ \frac{d\nu_{h_1}} {d\nu_{h_2}} \biggr] (\gamma, \sigma,
\beta_0, \beta_{\gamma})
&=& \biggl( \frac{w_1}{w_2}
\biggr)^{q_{\gamma}} \biggl( \frac{1 - w_1}{1 - w_2}
\biggr)^{q-q_{\gamma}} \nonumber\\[-8pt]\\[-8pt]
&&{}\times\frac{ \phi_{q_{\gamma}} (
\beta_{\gamma}; 0, g_1 \sigma^2 (X'_{\gamma} X_{\gamma})^{-1}
) } { \phi_{q_{\gamma}} ( \beta_{\gamma}; 0, g_2
\sigma^2 (X'_{\gamma} X_{\gamma})^{-1} ) },\nonumber
\end{eqnarray}
where $\phi_{q_{\gamma}}(u; a, V)$ is the density of the
$q_{\gamma}$-dimensional normal distribution with mean $a$ and
covariance $V$, evaluated at $u$ [\citet{Doss2007}]. It is
immediate that all formulas in Section \ref{secbf-pe} remain valid
if ratios of the form $\nu_h(\theta) / \nu_{h_1}(\theta)$ [see,
e.g., equation (\ref{ratio-est})] are replaced by the Radon--Nikodym
derivative $[ d\nu_h / d\nu_{h_1} ](\theta)$. Fortunately,
evaluation of (\ref{lr}) requires neither matrix inversion nor
calculation of a determinant, so can be done very quickly. Note
that in view of (\ref{lr}), it is not enough to have Markov chains
running on the $\gamma$'s and we need Markov chains running on the
$\theta$'s [or at least $(\gamma, \sigma, \beta_{\gamma})$].

There is a large literature on dealing with the hyperparameter in
models involving Zellner's $g$-prior [with or without the variable
inclusion line~(\ref{vsblm-c})]. Some of the proposals involve
putting a prior on $g$, or on both $g$ and $w$.
\citet{LiangEtal2008} propose and discuss priors on $g$; priors on
$w$ are generally taken to be beta distributions. Other proposals
give $g$ as a deterministic function of $m$ and $q$ [e.g., $g =
\max\{ m, q^2 \}$ in \citet{FernandezLeySteel2001}].
\citet{LiangEtal2008} contains an extensive and critical review of
the recommendations given in this literature. The most common
deterministic choice for $w$ is $w = 1/2$.
\citet{GeorgeFoster2000} recommend the empirical Bayes (EB)
approach for estimating the pair $(w, g)$: the marginal likelihood
of $(w, g)$ is computed over a grid, and the value of $(w, g)$ that
maximizes it is taken as the estimate of $(w, g)$. As with many
likelihood-based methods, special care needs to be taken when the
maximizing value is at the boundary. \citet{CuiGeorge2008} give
evidence that the EB method outperforms fully Bayes methods in this
problem. Unfortunately, the EB method is in general computationally
demanding because the likelihood is a sum over all $2^q$ models
$\gamma$, so it is practically feasible only for relatively small
values of $q$. Our methodology handles this problem by estimating
ratios of marginal likelihoods, that is, Bayes factors, and, besides
giving the maximizing values of $w$ and $g$, gives a plot which
shows the behavior of the Bayes factors for a wide range of other
values of $w$ and $g$.

We illustrate our methods on the US crime data of
\citet{Vandaele1978}, which can be found in the R library
\texttt{MASS} under the name \texttt{UScrime}. This data set seems
ideal, because it has been studied in several papers already, so we
can compare our results with previous analyses, and also because its
modest size enables a closed-form calculation of the marginal
likelihood $m_h$, so we can compare our estimates with the gold
standard. The data set gives, for each of $m = 47$ states of the
USA, the crime rate, defined as number of offenses per $100\mbox{,}000$
individuals (the response variable), and $q = 15$ predictors
measuring different characteristics of the population, such as
average number of years of schooling, average income, unemployment
rate, etc.

\begin{figure}

\includegraphics{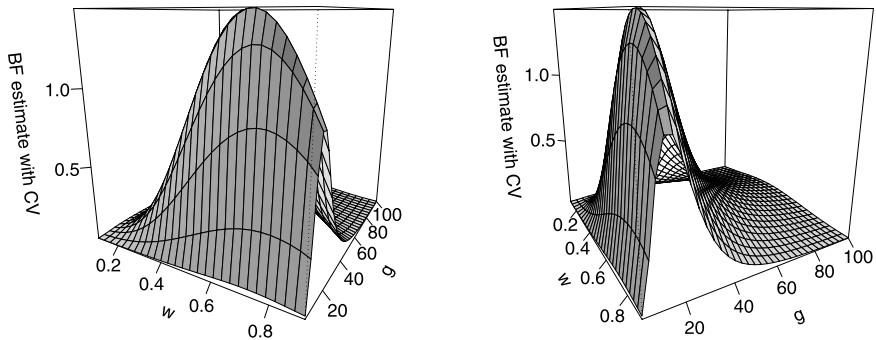}

\caption{Estimates of Bayes factors for the US crime data. The
plots give two different views of the graph of the Bayes factor as
a function of $w$ and $g$ when the baseline value of the
hyperparameter is given by $w = 0.5$ and $g = 15$. The estimate
is (\protect\ref{est-bf-complex-cv}), which uses control variates.}
\label{figbf}
\end{figure}

To be consistent with what is done in the literature, we applied a
log transformation to all variables, except the indicator variable.
We took the baseline hyperparameter to be $h_1 = (w_1, g_1) = (0.5,
15)$, and our goal was to estimate $B(h, h_1)$ for the $924$ values
of $h$ obtained when $w$ ranges from $0.1$ to $0.91$ by increments
of $0.03$, and $g$ ranges from $4$ to $100$ by increments of~$3$.
We used (\ref{est-bf-complex-cv}) and this estimate was based on
$16$ chains each of length $10\mbox{,}000$, corresponding to the skeleton
grid of hyperparameter values
%
%
\begin{equation}
\label{skeleton}
(w, g) \in\{ 0.3, 0.5, 0.6, 0.8 \} \times\{ 15, 50, 100, 225 \}
\end{equation}
for the stage $1$ samples, and $16$ new chains, each of length
$1\mbox{,}000$, corresponding to the same hyperparameter values, for the
stage $2$ samples. The plots in Figure \ref{figbf} give graphs of
the estimate (\ref{est-bf-complex-cv}) as $w$ and $g$ vary, from two
different angles. These indicate that values for $w$ around $0.65$
and for $g$ around $20$ seem appropriate, while values of $w$ less
than $0.3$ and values of $g$ greater than $60$ should be avoided. A
side calculation showed that, interestingly, for $g = \max\{ m, q^2
\}$ $(= 225)$, the estimate of $B( (w, g), (0.65, 20) )$ is
less than $0.008$ regardless of the value of $w$, so this choice
should not be used for this data set. With the long chains used and
the estimate that uses control variates, the Bayes factor estimates
in Figure \ref{figbf} are extremely accurate---root mean squared
errors are less than $0.04$ uniformly over the entire domain of the
plot and considerably less in the convex hull of the skeleton grid
(our calculation of the root mean squared errors used the
closed-form expression for the Bayes factors based on complete
enumeration). The figure took about a half hour to generate on an
Intel $2.8$ GHz Q$9550$ running Linux. (The accuracy we obtained is
overkill and the figure can be created in a few minutes if we use
more typical Markov chain lengths.)

%
%

Table \ref{tabcomparison} gives the posterior inclusion
probabilities for each of the fifteen predictors, that is, $P(\gamma_i
= 1 \mid y)$ for $i = 1, \ldots, 15$, under several models.
Line~$2$ gives the inclusion probabilities when we use
model (\ref{vsblm}) with the values $w = 0.65$ and $g = 20$, which
are the values at which the graph in Figure~\ref{figbf} attains its
maximum. Line $4$ gives the inclusion probabilities when the
hyper-$g$ prior ``HG$3$'' in \citet{LiangEtal2008} is used. As can
be seen, the inclusion probabilities we obtained under the EB model
are comparable to, but somewhat larger than, the probabilities when
the HG$3$ prior is used. This is not surprising since our model
allows $w$ to be chosen, and the data-driven choice gives a value
($0.65$) greater than the value $w = 0.5$ used in
\citet{LiangEtal2008}. [Table $2$ of \citet{LiangEtal2008} gives
a comparison of posterior inclusion probabilities for a total of ten
models taken from the literature.] Line $3$ of
Table \ref{tabcomparison} gives the inclusion probabilities under
model~(\ref{vsblm}) when we use $w = 0.5$ and the value of $g$ that
maximizes the likelihood with $w$ constrained to be $0.5$. It is
interesting to note that the inclusion probabilities are then
strikingly close to those under the HG$3$ model.

%
%
\begin{table}
\caption{Posterior inclusion probabilities for the fifteen
predictor variables in the US crime data set, under three models.
Names of the variables are as in Table $2$ of
Liang et al. (\protect\citeyear{LiangEtal2008}) (but all variables
except for the binary
variable S have been log transformed)}
\label{tabcomparison}
\begin{tabular*}{\tablewidth}{@{\extracolsep{\fill}}lcccccccc@{}}
\hline
& \textbf{Age} & \textbf{S} & \textbf{Ed} & \textbf{Ex0}
& \textbf{Ex1} & \textbf{LF} & \textbf{M} & \textbf{N}  \\
\hline
EB$(0.65,20)$ & $0.93$ & $0.39$ & $0.99$ & 0.$70$ & $0.51$ & $0.34$ &
$0.35$ & $0.52$
\\
EB$(0.5,20)$ & $0.85$ & $0.29$ & $0.97$ & 0.$67$ & $0.45$ & $0.22$ &
$0.22$ & $0.38$  \\
HG$3$ & $0.84$ & $0.29$ & $0.97$ & 0.$66$ & $0.47$ & $0.23$ & $0.23$ &
$0.39$
\end{tabular*}
\begin{tabular*}{\tablewidth}{@{\extracolsep{\fill}}lccccccc@{}}
\hline
& \textbf{NW} & \textbf{U1} & \textbf{U2} & \textbf{W} & \textbf{X}
& \textbf{Prison} & \textbf{Time} \\
\hline
EB$(0.65,20)$ & $0.83$ & $0.40$ & $0.76$ & $0.55$ & $1.00$ & $0.96$ & $0.55$\\
EB$(0.5,20)$ & $0.70$ & $0.27$ & $0.62$ & $0.38$ & $1.00$ & $0.90$ & $0.39$\\
HG$3$ & $0.69$ & $0.27$ & $0.61$ & $0.38$ & $0.99 $ & $0.89$ & $0.38$\\
\hline
\end{tabular*}
\vspace*{-4pt}
\end{table}

\citet{Buta2010} uses the estimates in Section \ref{secpe} to
produce plots of posterior inclusion probabilities for several of
the predictors, as $w$ and $g$ vary. The plots enable one to read
the posterior inclusion probabilities under various choices for $g$
and $w$ proposed in the literature, and also show that the extent to
which these probabilities change with the choices is striking.

%
%
\begin{figure}[b]
\vspace*{-4pt}
\includegraphics{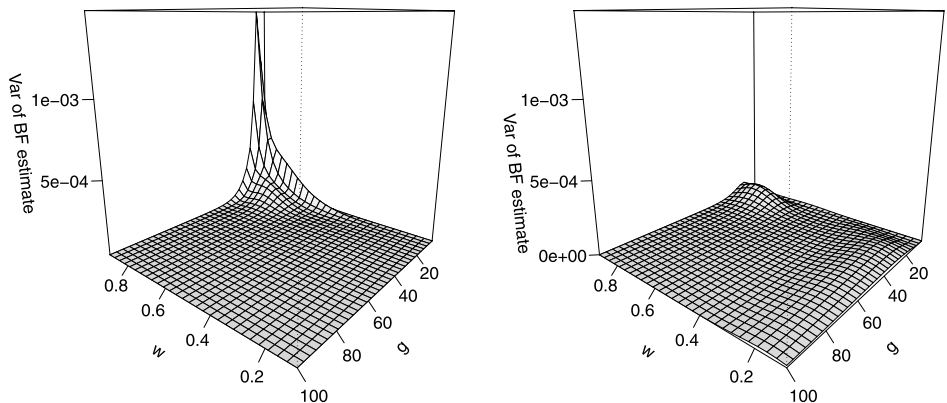}

\caption{Variance\vspace*{1pt} functions for two versions of
$\hat{I}_{{\hatbbeta(\hat{d})}}^{\hat{d}}$. The left panel is for
the estimate based on the skeleton (\protect\ref{skeleton}). The points
in this\vspace*{2pt} skeleton were shifted to better cover the problematic
region near the back of the plot ($g$ small and $w$ large),
creating the skeleton~(\protect\ref{skeleton-2}). The maximum
variance is
then reduced by a factor of $9$ (right panel).}
\label{figvar-est}
\end{figure}

Selection of the skeleton points was discussed at the end of
Section \ref{secest-var-sel-sp}, and we now return to this issue.
Consider the Bayes factor estimate based on the
skeleton~(\ref{skeleton}), which was chosen in an ad-hoc manner.
The left panel in Figure \ref{figvar-est} gives a plot of the
variance of this estimate, as a function of $h$. As can be seen
from the plot, the variance is greatest in the region where $g$ is
small and $w$ is large. We changed the skeleton
from (\ref{skeleton}) to
%
%
\begin{equation}
\label{skeleton-2}
(w, g) \in\{ 0.5, 0.7, 0.8, 0.9 \} \times\{ 10, 15, 50, 100 \}
\end{equation}
and reran the algorithm. The variance for the estimate based
on (\ref{skeleton-2}) is given by the right panel of
Figure \ref{figvar-est}, from which we see that the maximum
variance has been reduced by a factor of about $9$.\vadjust{\goodbreak}

\section{Discussion}
\label{secdisc}
The following fact is obvious, but it may be worthwhile to state it
explicitly. If $h_1$ is fixed, maximizing $B(h, h_1)$ and
maximizing the marginal likelihood $m_h$ are equivalent. Choosing
the value of $h$ that maximizes $m_h$ is by definition the empirical
Bayes method. Thus, the development in Section \ref{secbf-pe} can
be used to implement empirical Bayes methods.

Our methodology for dealing with the sensitivity analysis and model
selection problems discussed in Section \ref{secintro} can be
applied to many classes of Bayesian models. In addition to the
usual parametric models, we mention also Bayesian nonparametric
models involving mixtures of Dirichlet processes
[\citet{Antoniak1974}], in which one of the hyperparameters is the
so-called total mass parameter---very briefly, this hyperparameter
controls the extent to which the nonparametric model differs from a
purely parametric model. [Among the many papers that use such
models, we mention in particular \citet{BurrDoss2005}, who give a
more detailed discussion of the role of the total mass parameter.]
The approach developed in Sections \ref{secbf} and \ref{seccv} can
be used to select this parameter.

When the dimension of $h$ is low, it will be possible to plot $B(h,
h_1)$, or at least plot it as $h$ varies along some of its
dimensions. Empirical Bayes methods are notoriously difficult to
implement when the dimension of the hyperparameter $h$ is high. In
this case, it is possible to use the methods developed in
Sections \ref{secbf} and \ref{seccv} to enable approaches based on
stochastic search algorithms. These require the calculation of the
gradient $\partial B(h, h_1) / \partial h$. We note that the same
methodology used to estimate $B(h, h_1)$ can also be used to
estimate its gradient. For example, in (\ref{est-bf-complex}),
$\nu_h (\theta_i^{(l)})$ is simply replaced by $\partial\nu_h
(\theta_i^{(l)}) / \partial h$.

\begin{appendix}\label{app}
\section*{Appendix}

\begin{pf*}{Proof of Theorem \protect\ref{thmbf-d-est}}
We begin by writing
%
%
\begin{eqnarray}
\label{decomp}\quad
&&\sqrt{n} \bigl( \hat{B}(h, h_1, \hat{d}) - B(h, h_1) \bigr)
\nonumber\\[-8pt]\\[-8pt]
&&\qquad=
\sqrt{n} \bigl( \hat{B}(h, h_1, \hat{d}) - \hat{B}(h, h_1, d)
\bigr) + \sqrt{n} \bigl( \hat{B}(h, h_1, d) - B(h, h_1) \bigr).
\nonumber
\end{eqnarray}
The second term on the right-hand side of the equation in (\ref{decomp})
involves randomness coming only from the second stage of sampling.
This term was analyzed by \citet{Doss2010}, who showed that it is
asymptotically normal, with mean $0$ and variance $\tau^2(h)$. The
first term ostensibly involves randomness from both stage $1$ and
stage $2$ sampling. However, as will emerge from our proof, the
randomness from stage $2$ is of lower order, and effectively all the
randomness is from stage $1$. This randomness is nonnegligible.
We mention here the often-cited work of \citet{Geyer1994} (whose
nice results we use in the present paper). In the context of a
setup very similar to ours, his Theorem~$4$ states that using an
estimated $d$ and using the true $d$ results in the same asymptotic
variance. From our proof [refer also to the extension of our
Theorem \ref{thmbf-d-est} to the case of a simple sample given in
\citet{Buta2010}], we see that this statement is not correct.

To analyze the first term on the right-hand side of (\ref{decomp}),
define the function $F(u) = \hat{B}(h, h_1, u)$, where $u = (u_2,
\ldots, u_k)'$ is a real vector with $u_l > 0, l = 2, \ldots, k$.
Then, by the Taylor series expansion of $F$ about $d$, we get
%
%
\begin{eqnarray}
\label{tayl-exp}
&&\sqrt{n} \bigl( \hat{B}(h, h_1, \hat{d}) - \hat{B}(h, h_1, d)
\bigr) \nonumber\\
&&\qquad= \sqrt{n} \bigl( F(\hat{d}) - F(d) \bigr) \\
&&\qquad= \sqrt{n} \nabla F(d)' (\hat{d} - d) + \frac{\sqrt{n}}{2}
(\hat{d} - d)' \nabla^{2} F(d^{*}) (\hat{d} - d),
\nonumber
\end{eqnarray}
where $d^{*}$ is between $d$ and $\hat{d}$.

First, we show that the gradient $\nabla F(d) = (\partial F(d) /
\partial d_2, \ldots, \partial F(d) / \partial d_k)'$ converges
almost surely to a finite constant. Recall that $c(h)$ is defined
in~(\ref{c}). For $j = 2, \ldots, k$, the $(j-1)$th
component of $\nabla F(d)$ converges almost surely since, with the
SLLN assumed to hold for the Markov chains used, we have
\[
[\nabla{F(d)}]_{j-1} = \sum_{l=1}^k \frac{1}{n_l}
\sum_{i=1}^{n_l} \frac{ a_j a_l \nu_h (\theta_i^{(l)})
\nu_{h_j}(\theta_i^{(l)}) } { d_j^2 ( \sum_{s=1}^k a_s
\nu_{h_s} (\theta_i^{(l)}) / d_s )^2 } \cas[c(h)]_{j-1}.
\]

Next, we show that the random Hessian matrix $\nabla^{2} F(d^*)$ of
second-order derivatives of $F$ evaluated at $d^*$ is bounded in
probability. To this end, it suffices to show that each element of
this matrix, say, $[\nabla^{2} F(d^*)]_{t-1,j-1}$, where $t, j \in\{
2, \ldots, k \}$, is $O_p(1)$. Since $\vectornorm{d^* - d} \leq
\vectornorm{\hat{d} - d} \cpshort0$, it follows that $d^* \cpshort
d$.

Let $\varepsilon\in(0, \min(d_2, \ldots, d_k))$. Then we have
$\Pr(\vectornorm{d^* - d} \leq\varepsilon) \rightarrow1$. We now
show that, on the set $\{ \vectornorm{d^* - d} \leq\varepsilon\}$,
$\nabla^{2} F(d^*)$ is bounded in probability. Let
\[
\mathcal{I} = I(\vectornorm{d^* - d} \leq\varepsilon).
\]
For $t \neq j$, we have
%
%
\begin{eqnarray}
\label{Op1tneqj}
&&
| [ \nabla^2 F (d^*)
]_{t-1,j-1} | \cdot\mathcal{I} \nonumber\\
&&\qquad= \sum_{l=1}^k \frac{2}{n_l}
\sum_{i=1}^{n_l} \frac{ a_j a_l a_t \nu_h
(\theta_i^{(l)}) \nu_{h_j}(\theta_i^{(l)})
\nu_{h_t}(\theta_i^{(l)}) } { {d_j^*}^2 {d_t^*}^2
(\sum_{s=1}^k a_s \nu_{h_s} (\theta_i^{(l)}) /
d_s^*)^3 } \cdot\mathcal{I} \nonumber\\
&&\qquad\leq \sum_{l=1}^k \frac{2}{n_l} \sum_{i=1}^{n_l}
\frac{ a_j a_l a_t \nu_h (\theta_i^{(l)})
\nu_{h_j}(\theta_i^{(l)}) \nu_{h_t}(\theta_i^{(l)})
} { (d_j - \varepsilon)^2 (d_t - \varepsilon)^2 [
\sum_{s=1}^k a_s \nu_{h_s} (\theta_i^{(l)}) / (d_s +
\varepsilon) ]^3 } \nonumber\\
&&\qquad\cas \sum_{l=1}^k B(h, h_l) \int\biggl\{ \frac{ a_j
a_l a_t \nu_{h_j}(\theta) \nu_{h_t}(\theta)
\nu_{h_l}(\theta) } { [ \sum_{s=1}^k a_s
\nu_{h_s} (\theta) / (d_s + \varepsilon) ]^3 }
\biggr\} \nu_{h,y}(\theta) \,d\theta\nonumber\\[-8pt]\\[-8pt]
&&\qquad\quad\hspace*{7pt}{} \times\frac{2} { (d_j -
\varepsilon)^2 (d_t - \varepsilon)^2}. \nonumber
\end{eqnarray}
Note that the expression inside the braces in (\ref{Op1tneqj}) is
clearly bounded above by a constant, so expression (\ref{Op1tneqj})
is finite. Similarly, for $t = j$, we can show that $|
[\nabla^2 F(d^*) ]_{j-1,j-1} | \cdot\mathcal{I}$ is $O_p(1)$.
Since $\Pr(\vectornorm{d^* - d} \leq\varepsilon) \rightarrow1$, it
follows that~$\nabla^{2} F(d^*)$ is bounded in probability. Now, by
combining (\ref{decomp}) and (\ref{tayl-exp}), we obtain
\begin{eqnarray*}
&&\sqrt{n}\bigl( \hat{B}(h,
h_1, \hat{d}) - B(h, h_1) \bigr) \\
&&\qquad= \sqrt{\frac{n} {N}} \nabla
F(d)' \sqrt{N} (\hat{d} - d) \\
&&\qquad\quad{} + \frac{1}{2 \sqrt{N}} \sqrt{\frac{n} {N}} \bigl[
\sqrt{N} (\hat{d} - d) \bigr]' \nabla^{2} F(d^*) \bigl[
\sqrt{N} (\hat{d} - d) \bigr] \\
&&\qquad\quad{} + \sqrt{n} \bigl( \hat{B}(h, h_1, d) - B(h, h_1)
\bigr) \\
&&\qquad = \sqrt{q} c(h)' \sqrt{N} (\hat{d} - d) + \sqrt{n} \bigl(
\hat{B}(h, h_1, d) - B(h, h_1) \bigr) + o_p(1),
\end{eqnarray*}
where the last line follows from the fact that $\nabla F(d) \cas
{c(h)}$ established earlier, the assumptions of
Theorem \ref{thmbf-d-est} that $\sqrt{n/N} \rightarrow\sqrt{q}$
and that $\sqrt{N} (\hat{d}\,{-}\,d)$ converges in distribution [hence
is $O_p(1)$]. Because the two sampling stages [for estimating $d$
and $B(h, h_1)$] are assumed to be independent, using the assumption
that $\sqrt{N} (\hat{d} - d) \cd\mathcal{N}(0, \Sigma)$ in conjunction
with the result $\sqrt{n} ( \hat{B}(h,\allowbreak h_1, d) - B(h, h_1)
) \cd\mathcal{N}(0, \tau^2(h))$ established in Theorem $1$ of
\citet{Doss2010} under conditions \ref{A1} and \ref{A2}, we
conclude that
\[
\sqrt{n}\bigl( \hat{B}(h, h_1, \hat{d}) - B(h, h_1) \bigr) \cd
\mathcal{N}\bigl( 0, q c(h)' \Sigma c(h) + \tau^2(h) \bigr).
\]
\upqed
\end{pf*}
\begin{pf*}{Proof of Theorem \protect\ref{thmbf-d-est-cv}}
We begin by writing
%
%
\begin{equation}
\label{cov}\qquad
\sqrt{n} \bigl( \hat{I}_{{\hatbbeta(\hat{d})}}^{\hat{d}} - B(h,
h_1) \bigr) = \sqrt{n}\bigl(
\hat{I}_{{\hatbbeta(\hat{d})}}^{\hat{d}} - \hat{I}_{{\hatbbeta
(d)}}^d \bigr) + \sqrt{n} \bigl(\hat{I}_{{\hatbbeta(d)}}^d - B(h,
h_1) \bigr),
\end{equation}
where the second term on the right-hand side of (\ref{cov}) was analyzed
by \citet{Doss2010} who showed that it is asymptotically normal,
with mean $0$ and variance $\sigma^2(h)$. Our plan is to show that
$\hatbbeta(d)$ and $\hatbbeta(\hat{d})$ converge in probability to
the same limit, which we denote $\bbeta_{\lim}$. We then
expand the first term on the right-hand side of (\ref{cov}) by writing
%
%
\begin{eqnarray}
\label{first-term}
\sqrt{n} \bigl( \hat{I}_{{\hatbbeta(\hat{d})}}^{\hat{d}} -
\hat{I}_{\hatbbeta(d)}^d \bigr) &=& \sqrt{n} \bigl(
\hat{I}_{{\hatbbeta(\hat{d})}}^{\hat{d}} -
\hat{I}_{\bbeta_{\lim}}^{\hat{d}} \bigr) + \sqrt{n} (
\hat{I}_{\bbeta_{\lim}}^{\hat{d}} -
\hat{I}_{\bbeta_{\lim}}^{d} )\nonumber\\[-8pt]\\[-8pt]
&&{} + \sqrt{n} \bigl(
\hat{I}_{\bbeta_{\lim}}^{d} - \hat{I}_{\hatbbeta(d)}^{d}
\bigr).\nonumber
\end{eqnarray}
Our proof is organized as follows:
\begin{itemize}
\item We note that the third term on the right-hand side
of (\ref{first-term}) was shown to converge to $0$ in probability
by \citet{Doss2010}.
\item We will show that the first term on the right-hand side
of (\ref{first-term}) also converges to $0$ in probability.
\item The second term on the right-hand side of (\ref{first-term})
involves randomness from both stage $1$ and stage $2$. However,
we will show that the randomness from stage $2$ is asymptotically
negligible, and that this term is asymptotically equivalent to an
expression of the form $w(h)' (\hat{d} - d)$, where $w(h)$ is a~%
deterministic vector. This will show that the second term is
asymptotically normal.
\end{itemize}
Now we prove that the first term on the right-hand side
of (\ref{first-term}) is $o_p(1)$, and, to do this, we begin by
showing that $\hatbbeta(d)$ and $\hatbbeta(\hat{d})$ converge in
probability to the same limit. Let $\matbold{Z}$ be the $n
\times k$ matrix whose transpose is
%
%
\begin{equation}
\label{z-mat-t}
\matbold{Z}' =
\pmatrix{
1 & \cdots& 1 & 1 & \cdots& 1 & \cdots& 1 & \cdots& 1
\vspace*{2pt}\cr
Z_{1,1}^{(2)} & \cdots& Z_{n_1,1}^{(2)} & Z_{1,2}^{(2)} & \cdots&
Z_{n_2,2}^{(2)} & \cdots& Z_{1,k}^{(2)} & \cdots& Z_{n_k,k}^{(2)}
\cr
\vdots& \ddots& \vdots& \vdots& \ddots& \vdots& \ddots&
\vdots& \ddots& \vdots
\vspace*{2pt}\cr
Z_{1,1}^{(k)} & \cdots& Z_{n_1,1}^{(k)} & Z_{1,2}^{(k)} & \cdots&
Z_{n_2,2}^{(k)} & \cdots& Z_{1,k}^{(k)} & \cdots&
Z_{n_k,k}^{(k)}}\hspace*{-32pt}
\end{equation}
and let $\matbold{Y}$ be the vector
%
%
\begin{equation}
\label{Y-vect}
\matbold{Y} = ( Y_{1,1}, \ldots, Y_{n_1,1}, Y_{1,2},
\ldots, Y_{n_2,2}, \ldots, Y_{1,k}, \ldots, Y_{n_k,k} )'.
\end{equation}
Let $\hat{\matbold{Z}}$ be the $n \times k$ matrix corresponding
to $\matbold{Z}$ when we replace $d$ by $\hat{d}$. Similarly,
$\hat{\matbold{Y}}$ is like $\matbold{Y}$, but using $\hat{d}$
for $d$.

For fixed $j, j' \in\{ 2, \ldots, k \}$, consider the function
%
%
\begin{equation}
\label{G}
G(u) = \frac{1}{n} \sum_{l=1}^k \sum_{i=1}^{n_l} \frac{
\nu_{h_j}(\theta_i^{(l)}) / u_j - \nu_{h_1}(\theta_i^{(l)}) } {
\sum_{s=1}^k a_s \nu_{h_s}(\theta_i^{(l)}) / u_s } \cdot\frac{
\nu_{h_{j'}}(\theta_i^{(l)}) / u_{j'} - \nu_{h_1}(\theta_i^{(l)})
} { \sum_{s=1}^k a_s \nu_{h_s}(\theta_i^{(l)}) / u_s },\hspace*{-25pt}
\end{equation}
where $u = (u_2, \ldots, u_k)'$ and $u_l > 0$, for
$l = 2, \ldots, k$. [On the right-hand side of~(\ref{G}), $u_1$ is
taken to
be $1$.] Note that setting $u = d$ gives
\[
G(d) = \frac{1}{n} \sum_{l=1}^k \sum_{i=1}^{n_l} Z_{i,l}^{(j)}
Z_{i,l}^{(j')}.
\]
By the mean value theorem, we know that there exists a $d^*$ between
$d$ and~$\hat{d}$ such that
\[
G(\hat{d}) = G(d) + \nabla G(d^*)' (\hat{d} - d) =
\matbold{R}_{j,j'} + \nabla G(d^*)' (\hat{d} - d) + o_p(1).
\]
Note that\vspace*{1pt} the last equality above comes from applying the SLLN.
An argument similar to that used in Theorem \ref{thmbf-d-est} to
show that $\nabla^2 F(d^*) = O_p(1)$ can now be applied to show that
$\nabla G(d^*) = O_p(1)$.

Therefore,
\begin{eqnarray*}
G(\hat{d}) &=& \matbold{R}_{j,j'} + \nabla G(d^*)' (\hat{d} - d)
+ o_p(1) \\
&=& \matbold{R}_{j,j'} + O_p(1) o_p(1) + o_p(1) \cpshort
\matbold{R}_{j,j'}.
\end{eqnarray*}
Similar arguments extend to the case $j = 1$ or $j' = 1$. By the
fact that $\matbold{R}$ is assumed invertible, we have
%
%
\begin{equation}
\label{conv-inv-mat-hat}
n (\hat{\matbold{Z}}' \hat{\matbold{Z}})^{-1} \cp
\matbold{R}^{-1}.
\end{equation}
In a similar way, it can be shown that
%
%
\begin{equation}
\label{z-prime-y-hat}
\hat{\matbold{Z}}' \hat{\matbold{Y}} / n \cp\matbold{v},
\end{equation}
where $\matbold{v}$ is the same limit vector to which
$\matbold{Z}' \matbold{Y} / n$ has been proved to converge in
\citet{Doss2010}. Combining (\ref{conv-inv-mat-hat})
and (\ref{z-prime-y-hat}), we have
\[
( \hat{\beta}_0(\hat{d}), \hatbbeta(\hat{d}) ) = [
n (\hat{\matbold{Z}}' \hat{\matbold{Z}})^{-1} ] [
\hat{\matbold{Z}}' \hat{\matbold{Y}} / n ] \cp
(\beta_{0,\lim}, \bbeta_{\lim}) = \matbold{R}^{-1}
\matbold{v}.
\]
Let $e(j, l) = E( Z_{1,l}^{(j)} )$. We now have
%
%
\begin{eqnarray}
\label{newlabel}
\sqrt{n} \bigl( \hat{I}_{{\hatbbeta(\hat{d})}}^{\hat{d}} -
\hat{I}_{\bbeta_{\lim}}^{\hat{d}}
\bigr) & = & \sum_{j=2}^k \bigl( \beta_{j,\lim} -
\hat{\beta}_j(\hat{d}) \bigr) \Biggl( \sum_{l=1}^k a_l
n^{1/2} \sum_{i=1}^{n_l} \biggl( \frac{
\hat{Z}_{i,l}^{(j)} - e(j, l) } {n_l} \biggr) \Biggr)
\nonumber\hspace*{-30pt}\\[-8pt]\\[-8pt]
& = & \sum_{j=2}^k o_p(1) \Biggl( \sum_{l=1}^k a_l n^{1/2}
\sum_{i=1}^{n_l} \biggl( \frac{ \hat{Z}_{i,l}^{(j)} -
e(j, l) } {n_l} \biggr) \Biggr).\nonumber\hspace*{-30pt}
\end{eqnarray}
To show that (\ref{newlabel}) converges to $0$ in probability, it
suffices to show that for each $l$ and $j$
%
%
\begin{equation}
\label{cond}
n_l^{1/2} \sum_{i=1}^{n_l} \biggl( \frac{ \hat{Z}_{i,l}^{(j)} -
e(j, l) } {n_l} \biggr) = O_p(1).
\end{equation}
For fixed $j \in\{ 2, \ldots, k \}$ and $l \in\{ 1, \ldots, k \}$,
define
\[
H(u) = n_l^{-1/2} \sum_{i=1}^{n_l} \frac{ \nu_{h_j}
(\theta_i^{(l)}) / u_j - \nu_{h_1}(\theta_i^{(l)}) } {
\sum_{s=1}^k a_s \nu_{h_s}(\theta_i^{(l)}) / u_s }
\]
for $u = (u_2, \ldots, u_k)'$ with $u_l > 0, l = 2, \ldots, k$,
$u_1 = 1$. Note that $H(d) = n_l^{-1/2}\times \sum_{i=1}^{n_l}
Z_{i,l}^{(j)}$. To see why (\ref{cond}) is true, we begin by writing
%
%
\begin{eqnarray}
\label{tayl-H}
n_l^{1/2} \sum_{i=1}^{n_l} \biggl( \frac{ \hat{Z}_{i,l}^{(j)} -
e(j, l) } {n_l}
\biggr) & = & n_l^{1/2} \sum_{i=1}^{n_l} \biggl( \frac{
\hat{Z}_{i,l}^{(j)} - Z_{i,l}^{(j)} } {n_l} \biggr)\nonumber\\
&&{} +
n_l^{1/2}\sum_{i=1}^{n_l} \biggl( \frac{Z_{i,l}^{(j)}
- e(j, l)} {n_l} \biggr) \\
& = & H(\hat{d}) - H(d) + O_p(1).\nonumber
\end{eqnarray}
Note that the fact that $n_l^{1/2} \sum_{i=1}^{n_l} (
[Z_{i,l}^{(j)} - e(j, l)] / n_l ) = O_p(1)$, which was used to
establish the second equality in (\ref{tayl-H}), is proved in
\citet{Doss2010}. Now, applying the mean value theorem to the
function $H$, we know that there exists a point $d^*$ between $d$
and $\hat{d}$ such that (\ref{tayl-H}) becomes
%
%
\begin{eqnarray}
\label{tayl-H1}
n_l^{1/2} \sum_{i=1}^{n_l} \biggl( \frac{\hat{Z}_{i,l}^{(j)} -
e(j, l)} {n_l}
\biggr) & = & \nabla H(d^*)' (\hat{d} - d) + O_p(1) \nonumber\\
& = & \sqrt{a_l} \sqrt{\frac{n}{N}} n_l^{-1/2} \nabla
H(d^*)' \sqrt{N} (\hat{d} - d)\\
&&{} + O_p(1),\nonumber
\end{eqnarray}
so that the right-hand side of (\ref{tayl-H1}) is $O_p(1)$. We now
consider $\sqrt{n} ( \hat{I}_{\bbeta_{\lim}}^{\hat{d}} -
\hat{I}_{\bbeta_{\lim}}^{d} )$, the middle term
in (\ref{first-term}). Define
\[
K(u) = \frac{1}{n} \sum_{l=1}^k \sum_{i=1}^{n_l} \Biggl( \frac{
\nu_h(\theta_i^{(l)}) } { \sum_{s=1}^k a_s \nu_{h_s}
(\theta_i^{(l)}) / u_s } - \sum_{j=2}^k \beta_{j,\lim}
\frac{ \nu_{h_j}(\theta_i^{(l)}) / u_j - \nu_{h_1}(\theta_i^{(l)})
} { \sum_{s=1}^k a_s \nu_{h_s}(\theta_i^{(l)}) / u_s } \Biggr),
\]
where $u = (u_2, \ldots, u_k)'$, and $u_l > 0$ for $l = 2, \ldots,
k$. By the Taylor series expansion, we have
%
%
\begin{eqnarray}
\label{first-term-mid}
\sqrt{n} ( \hat{I}_{\bbeta_{\lim}}^{\hat{d}} -
\hat{I}_{\bbeta_{\lim}}^{d}
) &=& \sqrt{n} \nabla K(d)' (\hat{d} - d)\nonumber\\[-8pt]\\[-8pt]
&&{} + \sqrt{n}
\tfrac{1}{2} (\hat{d} - d)' \nabla^2 K(d^*)(\hat{d} -
d),\nonumber
\end{eqnarray}
where $d^*$ is between $\hat{d}$ and $d$. We now consider $\nabla
K(d)$. For $t = 2, \ldots, k$ we have
\[
[\nabla K(d)]_{t-1} \cas[w(h)]_{t-1},
\]
where $[w(h)]_{t-1}$ was defined in (\ref{grad-K-paper}). The
Hessian matrix $\nabla^2 K(d^*)$ can be shown to be bounded in
probability, using an argument similar to the one used in the proof
of Theorem \ref{thmbf-d-est}. Therefore, using the fact that
$\nabla^2 K(d^*)$ is bounded in probability, we can now
rewrite (\ref{first-term-mid}) as
\begin{eqnarray*}
\sqrt{n} ( \hat{I}_{\bbeta_{\lim}}^{\hat{d}}
- \hat{I}_{\bbeta_{\lim}}^{d}
) & = & \sqrt{\frac{n}{N}} w(h)' \sqrt{N} (\hat{d} - d) \\
&&{} + \sqrt{\frac{n}{N}}
\frac{1}{2\sqrt{N}} \sqrt{N}(\hat{d} - d)' O_p(1)
\sqrt{N} (\hat{d} - d) \\
& = &\sqrt{q} w(h)' \sqrt{N} (\hat{d} - d) + o_p(1).
\end{eqnarray*}
Together with (\ref{cov}), this gives
\begin{eqnarray*}
\sqrt{n} \bigl( \hat{I}_{{\hatbbeta(\hat{d})}}^{\hat{d}} - B(h,
h_1) \bigr) & = & \sqrt{q} w(h)' \sqrt{N} (\hat{d} - d) + \sqrt{n}
\bigl( \hat{I}_{{\hat{\bbeta}({d})}}^d
- B(h, h_1) \bigr) + o_p(1) \\
&\cd& \mathcal{N} \bigl( 0, qw(h)'
\Sigma w(h) + \sigma^2(h) \bigr)
\end{eqnarray*}
by the independence of the two stages of sampling, the assumption
that $\sqrt{N} (\hat{d} - d)$ is asymptotically normal with mean $0$
and variance $\Sigma$, and the result from \citet{Doss2010} that
$\sqrt{n} ( \hat{I}_{\hatbbeta(d)}^{d} - B(h, h_1) )$ is
asymptotically normal with mean $0$ and variance $\sigma^2(h)$.
\end{pf*}
\begin{pf*}{Proof of Theorem \protect\ref{thmpe-d-est}}
First, we note that
%
%
\begin{eqnarray}
\label{decomp-pe}
\sqrt{n} \bigl( \hat{I}^{[f]} (h, \hat{d}) - I^{[f]}(h) \bigr) &=&
\sqrt{n} \bigl( \hat{I}^{[f]}(h, \hat{d}) - \hat{I}^{[f]}(h, d)
\bigr) \nonumber\\[-8pt]\\[-8pt]
&&{}+ \sqrt{n} \bigl( \hat{I}^{[f]}(h, d) - I^{[f]}(h) \bigr).\nonumber
\end{eqnarray}
We begin by analyzing the second term on the right-hand side
of (\ref{decomp-pe}), which only involves randomness from the second
stage of sampling, and show that it is asymptotically normal. As
for the first term, a closer examination reveals that it is also
asymptotically normal, with all its randomness coming from
stage $1$. The asymptotic normality of the sum of these two terms
then follows immediately from the independence of the two stages of
sampling.

Note that $\sum_{l=1}^k a_l E( Y_{1,l}^{[f]} ) =
I^{[f]}(h) \cdot B(h, h_1)$, and, in particular, when $f \equiv1$,
this gives $\sum_{l=1}^k a_l E(Y_{1,l}) = B(h, h_1)$. Also, we have
%
%
\begin{eqnarray}
\label{biv-conv}
&&
n^{1/2}
\pmatrix{\displaystyle
\frac{1}{n} \sum_{l=1}^k \sum_{i=1}^{n_l} Y_{i,l}^{[f]} -
I^{[f]}(h) \cdot B(h, h_1) \vspace*{2pt}\cr
\displaystyle \frac{1}{n} \sum_{l=1}^k \sum_{i=1}^{n_l} Y_{i,l} - B(h, h_1)}
\nonumber\\
&&\qquad= n^{1/2}
\pmatrix{\displaystyle
\frac{1}{n} \sum_{l=1}^k \sum_{i=1}^{n_l} Y_{i,l}^{[f]} -
\sum_{l=1}^k a_l E\bigl( Y_{1,l}^{[f]} \bigr) \vspace*{2pt}\cr
\displaystyle \frac{1}{n} \sum_{l=1}^k \sum_{i=1}^{n_l} Y_{i,l} -
\sum_{l=1}^k a_l E(Y_{1,l})}
\\
&&\qquad= \sum_{l=1}^k {a_l}^{1/2} \cdot\frac{1}{{n_l}^{1/2}}
\sum_{i=1}^{n_l} \biggl[
\pmatrix{\displaystyle
Y_{i,l}^{[f]} \cr
Y_{i,l}}
-
\pmatrix{
E\bigl( Y_{1,l}^{[f]} \bigr) \cr
E(Y_{1,l})}
\biggr].\nonumber
\end{eqnarray}
By condition (\ref{mom-cond}), assumption \ref{A2} of
Theorem \ref{thmbf-d-est}, and the assumed geometric ergodicity and
independence of the $k$ Markov chains used, the vector
in (\ref{biv-conv}) converges in distribution to a normal random
vector with mean $0$ and covariance matrix $\Gamma(h)$ where
$\Gamma(h)$ is defined in (\ref{Gamma}). Since $\hat{I}^{[f]}(h,
d)$ is given by the ratio (\ref{Ihat-f-hd}), in view
of (\ref{biv-conv}), its asymptotic distribution may be obtained by
applying the delta method to the function $g(u, v) = u/v$. This
gives $\sqrt{n} ( \hat{I}^{[f]}(h, d) - I^{[f]}(h) ) \cd
\mathcal{N}(0, \rho(h))$, where $\rho(h)$ is given in (\ref{rho}).

We now consider the first term on the right-hand side
of (\ref{decomp-pe}). Define
\[
L(u) = \frac{ \sum_{l=1}^k \sum_{i=1}^{n_l}
({f(\theta_i^{(l)}) \nu_h (\theta_i^{(l)})}/ { \sum_{s=1}^k a_s
\nu_{h_s} (\theta_i^{(l)}) / u_{s} }) } {
\sum_{l=1}^k \sum_{i=1}^{n_l} ({\nu_h (\theta_i^{(l)})} /{
\sum_{s=1}^k a_s \nu_{h_s} (\theta_i^{(l)}) / u_{s} }) }
\]
for $u = (u_2, \ldots, u_k)'$ with $u_l > 0$ for $l = 2, \ldots,
k$. Then
\[
L(d) = \hat{I}^{[f]}(h, d) = \frac{ \sum_{l=1}^k \sum_{i=1}^{n_l}
Y_{i,l}^{[f]} } { \sum_{l=1}^k \sum_{i=1}^{n_l} Y_{i,l} }
\]
and $\sqrt{n}( \hat{I}^{[f]}(h, \hat{d}) - \hat{I}^{[f]}(h, d)
) = \sqrt{n} ( L(\hat{d}) - L(d) )$. Now, by the
Taylor series expansion of $L$ about $d$, we get
\[
\sqrt{n}\bigl( \hat{I}^{[f]}(h, \hat{d}) - \hat{I}^{[f]}(h, d)
\bigr) = \sqrt{n} \nabla L({d})' ({\hat{d}} - d) + \frac{\sqrt{n}}
{2} ({\hat{d}} - d)' \nabla^2 L(d^*) (\hat{d} - d),
\]
where $d^{*}$ is between ${d}$ and ${\hat{d}}$. First, we show that
the gradient $\nabla L({d})$ converges almost surely to a finite
constant vector by proving that each one of its components,
$[L(d)]_{j-1}, j = 2, \ldots, k$, converges almost surely. We
have
\[
[\nabla L(d)]_{j-1} \cas[v(h)]_{j-1},\qquad j = 2, \ldots, k,
\]
where $[v(h)]_{j-1}$ is given in (\ref{cas-grad-pe-paper}). As in
the proof of Theorem \ref{thmbf-d-est}, it can be shown that each
element of the second-derivative matrix $\nabla^{2} L(d^*)$ is
$O_p(1)$. Now, we can rewrite (\ref{decomp-pe}) as
\begin{eqnarray*}
&&\sqrt{n} \bigl( \hat{I}^{[f]}(h, \hat{d}) - I^{[f]}(h)
\bigr) \\
&&\qquad = \sqrt{\frac{n} {N}} \nabla L(d)' \sqrt{N} (\hat{d} - d)
+ \sqrt{n} \bigl( \hat{I}^{[f]}(h, d) - I^{[f]}(h)
\bigr) \\
&&\qquad\quad{} + \frac{1}{2 \sqrt{N}} \sqrt{\frac{n}
{N}} \bigl[ \sqrt{N} (\hat{d} - d) \bigr]' \nabla^{2}
L(d^*) \bigl[ \sqrt{N} (\hat{d} - d) \bigr] \\
&&\qquad = \sqrt{q} v(h)' \sqrt{N} (\hat{d} - d) +
\sqrt{n} \bigl( \hat{I}^{[f]} (h, d) - I^{[f]}(h)
\bigr) + o_p(1).
\end{eqnarray*}
Since the two sampling stages are assumed to be independent, we
conclude that
\[
\sqrt{n} \bigl( \hat{I}^{[f]}(h, \hat{d}) - I^{[f]}(h) \bigr) \cd
\mathcal{N} \bigl( 0, q v(h)' \Sigma v(h) + \rho(h) \bigr).
\]
\upqed
\end{pf*}
\end{appendix}

\section*{Acknowledgments}

We thank the reviewers for their careful reading and helpful comments.
We are especially grateful to the Associate Editor for a~very thorough
report and for suggestions which led to several improvements in the
paper.

\begin{supplement}[id=suppA]
\stitle{Additional technical details}
\slink[doi]{10.1214/11-AOS913SUPP} 
\sdatatype{.pdf}
\sfilename{aos913\_supp.pdf}
\sdescription{We show that when estimating the Bayes factors using
control variates, the estimate that is optimal when the
samples\vadjust{\goodbreak}
are i.i.d. sequences is no longer optimal when the samples are Markov
chains. We also give technical arguments regarding the
consistency of spectral estimates of the variance of our
estimators.}
\end{supplement}

%

\printaddresses

\end{document}